\documentclass[red,11pt,a4paper]{article}

\usepackage{cite}
\usepackage{comment}
\usepackage{amsmath}
\usepackage{amscd}
\usepackage{amssymb}
\usepackage[pdftex]{color,graphicx,hyperref}
\usepackage[utf8]{inputenc}
\usepackage{latexsym}
\usepackage{color}
\usepackage{cases}
\usepackage{xfrac}
\usepackage[normalem]{ulem}
\usepackage{environ}

\bf
\usepackage{ika}

\newcommand{\T}{\mathbb{T}}

\newcommand{\sqvr}{\sqrt{\varrho}}
\newcommand{\sqvrb}{\sqrt{\bar{\varrho}}}
\newcommand{\sqvrf}{\sqrt{\frac{\varrho}{\bar{\varrho}}}}
\newcommand{\sqvrft}{\sqrt{\frac{s}{\bar{\varrho}}}}
\newcommand{\Nu}{\mathcal{V}}
\newcommand{\Nutx}{\mathcal{V}_{t,x}}
\newcommand{\tvw}{{\tilde{\vw}}}
\newcommand{\lel}{\left\langle}
\newcommand{\ril}{\right \rangle}
\newcommand{\levert}{\left\vert}
\newcommand{\rivert}{\right \vert}
\newcommand{\levertl}{\left\Vert}
\newcommand{\rivertl}{\right \Vert}
\newcommand{\rtrw}{{\sqrt{s} \vc{v}}}
\newcommand{\DQ}{\mathbf{F} }
\newcommand{\intTAO}[1]{\int_0^\tau\!\! \int_{\mathbb{\T}^d} #1 \ \dxdt}

\NewEnviron{AS}{%
		\begin{align}
			\begin{split}
				\BODY
			\end{split}
		\end{align}
}
\title{Analysis of the generalised Aw-Rascle model}
\author{Nilasis Chaudhuri,$^*\;$
Piotr Gwiazda,$^\dagger\;$
Ewelina Zatorska$^\ddagger\;$
}

\date{\today}

\begin{document}
\maketitle

{
\footnotesize
\centerline{$^*\;$Department of Mathematics, Imperial College London}
\centerline{South Kensington Campus -- SW7 2AZ, London, UK}
\centerline{\small \texttt{n.chaudhuri@imperial.ac.uk}}

\bigbreak
\centerline{$^\dagger\;$Institute of Mathematics Polish Academy of Sciences}
\centerline{ul. \'Sniadeckich 8, 00-656 Warszawa, Poland}
\centerline{\small \texttt{pgwiazda@mimuw.edu.pl}}

\bigbreak
\centerline{$^\ddagger\;$Department of Mathematics, Imperial College London}
\centerline{South Kensington Campus -- SW7 2AZ, London, UK}
\centerline{\small \texttt{e.zatorska@imperial.ac.uk}}

}

\bigbreak

\begin{abstract}

We consider the multi-dimensional generalization of the Aw-Rascle system for vehicular traffic. For arbitrary large initial data and the periodic boundary conditions, we
prove the existence of global-in-time measure-valued solutions. We also show, using the relative energy technique, that the measure-valued solutions coincide with the classical solutions as long as the latter exist.

\end{abstract}




\bigskip

\maketitle
\section{Introduction}
The macroscopic modelling of traffic goes back to the work of Lighthill and  Whitham \cite{LW} and Richards \cite{R} from the late 50s.
Their description was using a single continuity equation for the conservation of mass
\eqh{\pt\vr+\px(\vr V(\vr))=0,}
in which the {\emph equilibrium} state velocity $V=V(\vr)$ was given explicitly. An additional equation for velocity, similar to that in the fluid models, appeared nearly 20 years later in the work of Payne \cite{Payne} and Whitham \cite{Whitham}. In 1995 Daganzo published a ``requiem" for the fluid models of traffic  \cite{Daganzo}, pointing out some of their serious drawbacks. They were subsequently remedied in the memorable work  of Aw and Rascle   \cite{AR} who ``resurrected"  the idea of using the second order fluid-like system to describe the traffic. The connection between the Aw-Rascle model and the microscopic {\emph{Follow-the-Leader}} was established in \cite{AwKlar} through a scaling limit.

The Aw-Rascle model in its homogeneous form (without relaxation term on the r.h.s.) describes the motion of cars on a single-lane road using two conservation laws:
\begin{subnumcases}{\label{AR}}
\pt \vr+\px(\vr u)=0,\label{AR1}\\
\pt(\vr w)+\px(\vr uw)=0,
\label{AR2}
\end{subnumcases} 
where $\vr(t,x)$ stands for the number of vehicles per unit length of road, a sort of one-dimensional density, and $u(t,x)$ stands for the velocity of the cars. The system \eqref{AR} resembles the compressible pressureless Euler equations, except that  the conservation of momentum involves not the actual velocity of motion $u$, but the preferred velocity $w$. Because the cars do not move freely on an empty road, these velocities  differ by a cost function (velocity offset) $p$ that depends on the concentration of the cars, i.e. $w=u+p(\vr)$. The use of the letter $p$ is not a coincidence. Here, $p(\vr)$ plays the same role as the pressure in the fluid mechanics: it propagates the flow perturbations. It should be noted, however, that because it is also the velocity offset, it has a different physical dimension. In one space dimension it does not strike as a problem, but it is not clear how to generalise the Aw-Rascle model to describe the multi-lane traffic.
One of possibilities\footnote{Private communications of P. Degond, A. Tosin, and  E. Zatorska} is to replace  the scalar cost function $p(\vr)$ by $\Grad p(\vr)$,
so that the velocity offset is also a vector, i.e. $\vw=\vu+\Grad p(\vr)$. We then arrive at the following multi-dimensional version of \eqref{AR}:
\begin{subnumcases}{\label{AR_multi}}
\pt \vr+\Div (\vr \vu)=0,\label{AR1_multi}\\
\pt(\vr \vw)+\Div(\vr \vu\otimes\vw)=0.
\label{AR2_multi}
\end{subnumcases} 
The first main result  of this work is the existence of global-in-time solutions to system \eqref{AR_multi} for arbitrary large initial data. We look for solutions in terms of $\vr$ and $\vw$ treating $\vu$ as given function, defined by
\eq{\label{def:w}
\vu=\vw-\Grad p(\vr), \quad \mbox{where}\quad p(\vr)=\vr^\gamma,\quad \gamma\geq1.
}
For simplicity, we restrict the analysis to periodic spacial domain, i.e. $\Omega=\mathbb{T}^d$, $d=2,3$. The proof of existence of solutions to the system \eqref{AR} consists of a multi-layer construction procedure leading eventually to the generalized measure-valued solution. 
A remarkable property of the measure-valued solutions is their usability in the construction of numerical algorithms. For example, in \cite{FMT} it was shown that the widely accepted entropy solutions to Euler system are not the ones observed as a limit of numerical schemes. More recently, the weak-strong uniqueness of dissipative measure-valued solutions was used to prove the convergence of finite volume numerical schemes for the Euler and the Navier-Stokes equations \cite{FL, FLM}. For an overview of results in this direction we refer to the recent monograph \cite{FLMS}.

Motivated by these findings, in our second main result we show that the measure-valued solutions to system \eqref{AR_multi} coincide with the classical solutions emanating from the same initial data, as long as the latter exist. This is the so-called weak(measure-valued)-strong uniqueness principle. Its proof relies on the application of the relative energy method. Roughly speaking, we show that the functional measuring the distance between the solutions can be controlled by the difference of the corresponding initial data. For two smooth solutions $ (\bvr, \bar{\vw}) $ and $ (\vr,\vc{w}) $ this functional has a form:
\eq{\label{def:rel_ent}
\mathcal{E}(\varrho, \vc{w} \mid \bvr , \bar{\vw}) : = \intO{\left(\frac{1}{2} \varrho \vert \vw- \bar{\vw} \vert ^2 +  E(\varrho)-E(\bvr)-E^\prime (\bvr) (\varrho-\bvr) \right)}
}
where
\eq{\label{def:E}
E(\vr)=\int_0^\vr p(s)\ {\rm d}s.}
{One of the main challenges in our proof is} to verify that one of the smooth solutions in this functional can be replaced by the measure-valued solution constructed in the first part.

The concept of relative energy was first introduced by Dafermos \cite{D1979} in the context of the hyperbolic conservation laws. It was later on adapted to prove weak-strong uniqueness (admissibility) of the weak solutions to the compressible Navier-Stokes equations by Feireisl, Jin, Novotn\'y \cite{FJN2012}, and by Feireisl, Novotn\'y and Sun in \cite{FNS2011}. We also refer to the review paper by Wiedemann \cite{W2018} for a  survey of recent results on weak-strong uniqueness for compressible and incompressible Euler and Navier-Stokes equations.
The weak-strong property of measure-valued solutions for the incompressible Euler system was first observed by Brennier et al. \cite{Brenier}. After that, analogous results were proven for the whole variety of fluid models \cite{GSW, ChTz,CDGS} including ultimately the compressible Navier-Stokes equations \cite{FGSW}, and general conservation laws \cite{GKS}.

One more observation that should be made at this stage is that system  \eqref{AR} resembles compressible pressureless Navier-Stokes equations with variable viscosity coefficient. Indeed, rewriting the equations in terms of $(\vr,u)$, we obtain:
\begin{subnumcases}{\label{NS}}
\pt \vr+\px(\vr u)=0,\label{NS1}\\
\pt(\vr u)+\px(\vr u^2)-\px(\vr^2p'(\vr)\px u)=0.
\label{NS2}
\end{subnumcases} 
Note that the viscosity coefficient $\vr^2p'(\vr)$ vanishes when density is equal to $0$. For that reason, system \eqref{NS} is sometimes called the degenerate Navier-Stokes system. The weak solutions to \eqref{NS} were obtained in the vanishing pressure limit by Haspot and Zatorska in \cite{HZ}. They used the artificial velocity reformulation \eqref{AR} to compare the profile of the density with the Barenblatt solutions to the porous medium equation. An analogue of this result in multi-dimensional case is due to Haspot \cite{PAM}. The two-velocity  structure of the degenerate compressible Navier-Stokes equations was used earlier by Bresch and Desjardins \cite{BD, BrDe} to derive additional estimate of the gradient of the density. This extra estimate is nowadays called the BD entropy, and the algebraic relation for the viscosity coefficients necessary to close this estimate is  called the BD relation. Later on, this relation was generalised, and the concept of $\kappa$-entropy two-velocity solutions was introduced in the series of works \cite{BrGiZa,BrDeZa}. The relative-entropy functional based on this concept was proposed for the degenerate compressible Navier-Stokes system  in \cite{BrNoVi2017}. The complete approximation scheme allowing to construct a weak solution in case of linear viscosity is due to Vasseur and Yu \cite{VaYu}. In all these works the additional estimate for the density compensates for the lack of information on the velocity vector field that would be independent of $\vr$. In a sense, the velocity vector field cases to be defined in the regions of vacuum. In a sense the momentum and the continuity equations play the opposite roles in comparison to compressible Navier-Stokes  equations with constant viscosity coefficients. A similar shift of regularity between the momentum and the continuity equations has been recently observed for the one-dimensional collective motion models, see for example \cite{CaWrZa2,kis,tad4}, {and in the so called Brenner model \cite{B1,B2,FV}.}

{\section{Definition of solution and main results}}
In this section, we first explain the reasons for working in the framework of generalized  measure-valued solution. After this, we introduce some basics of the Young measures theory, and explain the extensions necessary to define the measure-valued solutions to system  \eqref{AR_multi}.

\subsection{Motivation: a-priori estimates}
Below we derive some basic a-priori estimates for smooth solutions of system \eqref{AR_multi} supplemented with the initial data $\vr(0,x)=\vr_0$ and $(\vr\vw)(0,x)=\vc{m}_0$ such that
\eq{\label{MV:IE:2}
        &0\leq  \varrho_0\quad \text{ in } \T^d\quad \text{ and }\quad E_0=\intO{\left( \frac{1}{2} \frac{\vert \vc{m}_0 \vert^2}{\vr_0} + E(\vr_0)\right)} < \infty .
}
{The above choice of initial data includes  a possibility of vacuum (i.e. $ \vr_0=0$) in some subset of $ \T^d$, but the condition \eqref{MV:IE:2} ensures that $ \vc{m}_0=0$ on $\{ x\in \T^d| \vr_0(x)=0\}$.}
\begin{lemma}\label{Lemma:1}
Let $(\vr,\vw)$ be a smooth solution to \eqref{AR_multi} and let $\vu$ be given by \eqref{def:w}. Let the initial data satisfy \eqref{MV:IE:2}, then we have for any $p\in[2,\infty)$
\eq{\label{est_mas}
\|\vr\|_{L^\infty(0,T; L^1(\T^d)}\leq C,
}
\eq{\label{est_mom}
\|\vr|\vw|^p\|_{L^\infty(0,T; L^1(\T^d))}\leq C,
}
where the constant $C$ depends only on the initial data.
\end{lemma}
\pf The proof is quite straightforward. To obtain \eqref{est_mas} we integrate equation \eqref{AR1_multi} over space, and use the periodicity of the domain. To obtain \eqref{est_mom}, we multiply \eqref{AR2_multi} by $\vw|\vw|^{p-2}$ and \eqref{AR1_multi} by $\frac1p|\vw|^p$, we integrate over space and sum up the resulting expressions. $\Box$

\bigskip
There is another interesting way of thinking about system \eqref{AR_multi}, namely we can formally rewrite it in terms of $\vr$ and $\vw$ as:
\begin{subnumcases}{\label{PM}}
\pt \vr+\Div(\vr \vw)-\Div(\vr\Grad p)=0,\label{PM1}\\
\pt(\vr \vw)+\Div(\vr \vw\otimes\vw)=\Div(\vr \Grad p\otimes\vw).
\label{PM2}
\end{subnumcases} 
System \eqref{PM} is a porous medium equation coupled to convection equation with a drift term. This reformulation allows us to get further density estimates.

\begin{lemma}\label{Lemma:2}
Let assumptions of Lemma \ref{Lemma:1} be satisfied.
Then we have the following estimates
\begin{align}
   & \|E(\vr)\|_{L^\infty(0,T; L^1(\T^d))}\leq C, \label{eng:1}\\
   & \|\vr|\Grad p(\vr)|^2\|_{L^1((0,T)\times\T^d)}\leq C, \label{eng:2}
\end{align}
where the constant $C$ depends only on the initial data.
\end{lemma}
\pf
We multiply the continuity equation \eqref{PM1} by $p(\vr)$ and  integrate over space
\eqh{
\intO{\pt\vr\ p(\vr)}-\intO{\vr\vw\cdot\Grad p(\vr)}+\intO{\vr|\Grad p(\vr)|^2}=0.
}
Definition \eqref{def:E} implies that $E'(\vr)=p(\vr)$, and thus
\begin{align}\label{energy-eqn-2}
    \Dt\intO{E(\vr)}-\intO{\vr\vw\cdot\Grad p(\vr)}+\intO{\vr|\Grad p(\vr)|^2}=0.
\end{align}
Using Young's inequality with  $\delta<<1$ to estimate the second term, we obtain
\eqh{
\Dt\intO{E(\vr)}+\intO{\vr|\Grad p(\vr)|^2}\leq C_\delta\intO{\vr|\vw|^2}+\delta\intO{\vr|\Grad p(\vr)|^2}.
}
The last term can be thus absorbed by the left hand side (l.h.s), and we conclude the proof by integrating over time and using \eqref{est_mom} to bound the first term on the r.h.s. $\Box$

As a consequence of the above two lemmas, it is easy to check that all the terms of the system \eqref{AR_multi} are bounded in the sense of distributions. The $L^1$ bounds for all terms in the weak formulation of the system \eqref{AR_multi} could allow us to consider it as possible framework for the existence theory. However, the Lemmas \ref{Lemma:1} and \ref{Lemma:2} do not provide a-priori bound for $ \vw $ that would be independent of $\vr$. Indeed, the $L^p$ bound for $\sqrt{\vr}\vw$ is the maximum we can expect.  This causes that the nonlinear terms in the momentum equation \eqref{PM2} are difficult to handle, similarly as in the case of the compressible Euler system. To bypass this difficulty we employ the notion of Young measures which are tailored to identify the limits in the nonlinear terms. 

Below we recall some of the basic facts from the theory of Young measures.

{\subsection{Basic facts about the Young measures}\label{Ym}}

We first introduce some notation:
\begin{itemize}
    \item $ C_0(\R^d) $ denotes the closure
under the supremum norm of compactly supported, continuous functions on $ \R^d $,
that is the set of continuous functions on $ \R^d $ vanishing at infinity;
\item $\mathcal{M}(\R^d) $ denotes the dual space of $ C_0(\R^d) $ consisting of signed Radon measures with finite mass, equipped with the dual norm of total variation;
\item $\mathcal{M}^{+}(\R^d)$ denotes the cone of non-negative Radon measures on $\R^d$ and $\mathcal{P}(\R^d)$ indicates the space of probability measures, i.e. for $\nu \in \mathcal{P}(\R^d)\subset {\mathcal{M}^+} (\R^d)$ we have $\nu[\R^d]=1$. 
\end{itemize}

Let $ \Omega = (0,T)\times \T^d $ for  any $ T>0 $, and let $\{\vc{U}_n\}_{n\in\mathbb{N}}$ be a sequence of measurable functions 
$$  \vc{U}_n \colon  \Omega  \rightarrow \R^k .$$ 
With a slight abuse of notation we identify $\vc{U}_n$ with the probability measures  $$ \vc{U}_n \colon y(\in  \Omega ) \mapsto \delta_{\vc{U}_n(y)} \in \mathcal{P}(\R^k) \subset \mathcal{M}(\R^k). $$ 
The new sequence converges to a measure $\Nu$ parameterized by $y$ in the natural weak topology  
\[ \vc{U}_n \rightarrow  \Nu \text{ weak-(*)ly in } L^{\infty}_{\text{weak-(*)}} \big( \Omega  ;\mathcal{M}(\R^k )\big)=[L^1\big( \Omega ;C_0(\R^k)\big)]^{*} .\] 
Assuming that $ \Vert \vc{U}_n \Vert_{L^1( \Omega ;\R^k)} \leq C $ uniformly in $ n $, allows us to deduce that $\Nu$ is a probability measure, i.e.,
\[ \Nu \in  L^{\infty}_{\text{weak-(*)}} \big(\Omega ;\mathcal{P}(\R^k )\big), \] 
 called the {\emph{Young measure}} associated to the sequence  $ \{ \vc{U}_n \}_{n\in\mathbb{N}}$.

We now provide a few more important results on the Young measures theory. For a more complete overview we refer the reader, for example, to Pedregal \cite[Chapter 6]{P}, Malek et al. \cite[Chapter 3]{MNRR}, or a recent monograph by Feireisl et al. \cite[Chapter 5]{FLMS}.

If $\vc{U}_n \in L^p(\Omega;\R^k) \text{ for } \Vert \vc{U}_n \Vert_{L^p(\Omega;\R^k)}  \leq C$ uniformly in $ n$ and $1<p\leq \infty$, then for $1<p<\infty$, the weak limit of $\vc{U}_n $ (weak-(*), for $p=\infty$ ), denoted by $\vc{U}$, is given by the barycenter of the Young measure, i.e., 
 \[ \vc{U}(y)= \langle \Nu_{y}; {\vc{Id}} \rangle:=
 \int_{\R^{k}} \lambda \text{d}\Nu_{y}(\lambda). \]

 The situation becomes much more interesting when $\vc{U}_n \in L^1(\Omega;\R^k) \text{ and } \Vert \vc{U}_n \Vert_{L^1(\Omega;\R^k)}  \leq C$ uniformly in $n$. Following  Ball and Murat \cite[Section 3]{BM}, we conclude that  $\langle \Nu_{y}; {\vc{Id}} \rangle\left(=
 \int_{\R^{k}} \lambda \text{d}\Nu_{y}(\lambda)\right)$ coincides with the \textit{biting limit} of $ \vc{U}_n $, denoted by $ \vc{U} $, defined as follows. 
 There exists a function $\vc{U} \in L^1(\Omega;\R^k)$, a subsequence  $ \vc{U}_{n_k} $ of  $ \vc{U}_n $ and a non-increasing sequence of Borel subsets of $\Omega$, denoted by $\{E_l\}_{l\in \mathbb{N}}$ with $\lim_{l\rightarrow \infty} \vert E_l\vert =0$, such that 
 \[ \vc{U}_{n_k} \rightarrow \vc{U} \text{ weakly in }L^1(\Omega\setminus E_l;\R^k),\]
 as $n_k \rightarrow \infty $ for fixed $l$. More details on the biting limit and its identification are provided in \cite[Section 1]{BM}.

 Let us now understand how the Young measures can be applied to identify the limits of nonlinear functions of weakly converging sequences.\\
 For any continuous function  $ \vc{b}: \R^k \rightarrow \R^d $, with $ \Vert \vc{U}_n \Vert_{L^p(\Omega;\R^d)}+ \Vert \vc{b}(\vc{U}_n) \Vert_{L^p(\Omega;\R^d)}  \leq C \text{ and } p\geq 1 $ where $ C $ is independent of $ n $. Then the mapping
 \[ y \mapsto
 \int_{\R^{k}} \vc{b}(\lambda) \text{d}\Nu_{y}(\lambda) \] is well defined and finite for a.e. $y\in \Omega$. Here $ \{ \Nu_y\}_{y\in Q} $ is a Young measure generated by $\{ \vc{U}_n \}_{n\in \mathbb{N}} $ and we have
	\begin{align*}
		 y \mapsto \int_{\R^{k}} \vc{b}(\lambda) \text{d}\Nu_{y}(\lambda)  \in L^p(\Omega),
	\end{align*} see \cite[Proposition 5.2]{FLMS}. It is convenient to introduce the notation \begin{align}\label{Non-ym}
    \langle \Nu_{y}; \vc{b} \rangle := \int_{\R^{k}} \vc{b}(\lambda) \text{d}\Nu_{y}(\lambda)
\end{align}
Let us now consider the case that
	\begin{align*}
		\vc{U}_n \in L^1(\Omega;\R^k),\quad  \Vert \vc{U}_n \Vert_{L^1(\Omega;\R^k)}  \leq C,\quad \Vert b(\vc{U}_n ) \Vert_{L^1(\Omega;\R^d)} \leq C 
	\end{align*}
	uniformly w.r.t $n$. From the continuous embedding of $L^1(\Omega)$ into   $ \mathcal{M}(\Omega) $, we have $$ \vc{b}(\vc{U}_n) \rightarrow \Ov{\vc{b}(\vc{U})} \text{ weak-(*)ly in\ }  \mathcal{M}(\Omega;\R^k), $$ 
    which does not necessarily coincide with $\left\langle \Nu_y; \vc{b} \right\rangle$. Therefore, we define a \textit{Concentration Defect} as {the difference} 
    $$ \mathcal{R}^{\text{cd}}=  \Ov{\vc{b}(\vc{U})} -	\left\langle \Nu_y; \vc{b}\right\rangle \in \mathcal{M}(\Omega;\R^k) .$$
It is important to note that, if we assume
\begin{align*}
		\vc{U}_n \in L^p(\Omega;\R^k),\quad  \Vert \vc{U}_n \Vert_{L^p(\Omega;\R^k)}  \leq C,\quad \Vert b(\vc{U}_n ) \Vert_{L^p(\Omega;\R^k)} \leq C 
	\end{align*}
	uniformly w.r.t $n$ and $p>1$, we have $\mathcal{R}^{\text{cd}} =0 $.
 
\subsection{{Definition of the measure-valued  solution of \eqref{PM}}}
In this article, we will focus on the existence of solutions to the Aw-Rascle system \eqref{PM}, which is equivalent to \eqref{AR_multi} on the level of classical solutions.
Perhaps, the most natural idea would be to consider a measure-valued solution defined by a Young measure $ \Nu$ generated by the sequence $ \{\vr_n, \vr_n \vw_n\}_{n\in \mathbb{N}}$, in full analogy with the compressible Euler system, \cite[Chapter 2]{Basaric}.\par 
Unfortunately, in case of system \eqref{PM} the nonlinear term $\vrn \Grad p(\vrn)\otimes\vw_n $ is only a-priori bounded in $L^1$ and its limit cannot be identified in this space. This obstacle requires to consider a generalized Young measure involving the density $\vr$, the gradient of a certain function of the density 
\begin{align}\label{Q-defn}
	\Grad Q(\vr)=\sqrt{\vr}\Grad p(\vr), \text{ where } 	Q^\prime (\vr)= \sqrt{ \vr} p^\prime (\vr)\quad \text{ for }\quad \vr\geq 0,
\end{align}
and the ``weighted'' velocity  $\sqrt{\vr} \vw$.
Similar idea appeared in the articles of B\v rezina, Feireisl, and Novoton\' y \cite{BFN} and also in \cite{Ch20} devoted to the compressible Navier-Stokes system, but in the context of velocity and velocity gradient.\par

Summarising the above considerations, the solution to \eqref{AR_multi} will be the Young measure $ \Nu $, generated by $ \left\{(\vr_n, \sqrt{\vr_n} \vw_n, \Grad Q(\vr_n))\right\}_{n\in \mathbb{N}} $, where $(\vr_n, \vw_n)$ is the solution of a suitable approximation of \eqref{PM}.
Note however, that thanks to \eqref{PM1} and the uniform bounds from Lemma \ref{Lemma:2}, we gain some control over the oscillations of the sequences which approximate the density. In fact, we are able to show that these sequences converge strongly. This allows us to deduce that the measure associated to the nonlinear functions of the density can be identified with a Dirac delta. For further details, we refer the reader to Lemma \ref{Lemma:strong} below. \par 
{

Suppose that $ \left\{(\vr_n, \sqrt{\vr_n} \vw_n, \Grad Q(\vr_n))\right\}_{n\in \mathbb{N}} $ generates the Young measure $ \Nu $, i.e. 
\begin{align*}
		\Nu \in L^{\infty}_{\text{weak-(*)}} \big( (0,T)\times \T^d ;\mathcal{P}(\mathcal{F} )\big) \text{ where } \mathcal{F}= \left\{\vc{\lambda}  \mid \vc{\lambda} \in [0,\infty)\times  \R^d \times  \R^d \right\}.
	\end{align*}
	Consider a continuous function $ \vc{b}=\vc{b}(\vr_n, \sqrt{\vr_n} \vw_n, \Grad Q(\vr_n)) : \mathcal{F} \rightarrow \R^k$, which certain uniform $L^p$-bound with $p\geq 1$. In agreement with \eqref{Non-ym}, we  have \begin{align*}
    \lel \Nu_{t,x}; \vc{b} \ril = \int_{\mathcal{F}} \vc{b}(\vc{\lambda})  \;  \text{d} \Nu_{t,x}(\vc{\lambda}) , \quad \mbox{where}\quad \vc{\lambda} = (s, \sqrt{s}\vc{v}, \vc{F}).
\end{align*}
Now we want to consider a particular nonlinearity that appears in the convective term of \eqref{PM2}, that is $\vr \vw \otimes \vw $. In this case our
our nonlinear function can be written as
\[ \vc{b}(s, \sqrt{s}\vc{v}, \vc{F}) =\sqrt{s}\vc{v} \otimes \sqrt{s}\vc{v} . \] 
In what follows, we use a slight abuse of notation and simply  write
\begin{align*}
    \lel \Nu_{t,x}; \vc{b} \ril = \int_{\mathcal{F}} \sqrt{s}\vc{v} \otimes \sqrt{s}\vc{v} \;  \text{d} \Nu_{t,x}(s, \sqrt{s}\vc{v}, \vc{F}) := \lel \Nu_{t,x}; \rtrw \otimes  \rtrw  \ril .
\end{align*}
 Similarly, we can express the other nonlinear term $ \lel \Nu_{t,x}; \rtrw \otimes  \DQ  \ril $ in \eqref{PM2}. In the literature, the variables $(s, \sqrt{s}\vc{v}, \vc{F})$ are often called ``dummy'' variables. Moreover, in the above notation  the barycenter of the Young measure as \begin{align*}
    \lel \Nu_{t,x}; \vc{Id} \ril &= \int_{\mathcal{F}} \left( (s, \sqrt{s}\vc{v}, \vc{F}) \right) \text{d} \Nu_{t,x}(s, \sqrt{s}\vc{v}, \vc{F}) := \left( \lel \Nu_{t,x}; s  \ril ,\lel \Nu_{t,x};  \sqrt{s}\vc{v} \ril,\lel \Nu_{t,x};  \vc{F} \ril\right)  .
\end{align*} }

Let us also rewrite the assumptions on the initial data \eqref{MV:IE:2} in terms of our new unknowns  $(\vr,\sqrt{\vr} \vw, \Grad Q(\vr))$ generating the Young measure. From the estimate \eqref{MV:IE:2} it is evident that $ \vc{m}_0=0 $ a.e. on the set $\{ x\in \T^d \mid \vr_0(x) =0  \} $. 
Therefore, considering $\vw_0 = \frac{\vc{m}_0}{\vr_0} $, our bounds on initial data can be rewritten in the form:
\eq{\label{MV:IE:11}
        &0\leq  \varrho_0 \text{ in } \T^d\quad \text{ and }\quad E_0=\intO{\left( \frac{1}{2}\vr_0\vert \vw_0 \vert^2 + E(\vr_0)\right)} < \infty .
}

For the existence of solutions, one can consider even more general assumption  in terms of the Young measures, i.e.,  to assume for some $\Nu_{0} \in L^{\infty}_{\text{weak-(*)}} \big( \T^d ;\mathcal{P}([0,\infty )\times \R^d)\big) $, that
\begin{align}\label{MV:IE:1}
     0\leq \lel  \Nu_{0,x}; s \ril, \text{ for a.e.}~ x \in \T^d \text{ and }  \int_{ \T^d} \lel \Nu_{0,x}; \frac{1}{2}  \levert \rtrw \rivert^2 + E(s) \ril \dx < \infty ,
\end{align}
		where $ (s, \rtrw) $ are the dummy variables. Note that since the initial energy depends  only on $ \vr $ and $ \sqrt{\vr}\vw $, it is enough to consider the phase space of $ \Nu_{0,x} $ as $([0,\infty )\times \R^d) $. 

Our ultimate goal is to prove the weak(measure-valued)--strong uniqueness, and so, the initial data should coincide for both measure-valued and strong solutions.  Therefore, for this part we will use  \eqref{MV:IE:2}. In other words, we will assume \eqref{MV:IE:1} with the initial measure  $\Nu_{0,x} = \delta_{\{\vr_0(x), \sqrt{\vr_0}\vw_0(x)\}} \;\text{ for a.e. }x \in \T^d $.\par

Let us now provide a precise definition of the measure-valued solution to the system \eqref{PM}.

\begin{df}\label{MVdef}
		We say that a parametrized measure $\{ \Nu_{t,x} \}_{(t,x)\in (0,T)\times \Omega}$,
	\begin{align*}
		\Nu \in L^{\infty}_{\text{weak-(*)}} \big( (0,T)\times \T^d ;\mathcal{P}(\mathcal{F} )\big),
	\end{align*}
	on the phase space
	\begin{align*}
		\mathcal{F}= \left\{\left(s, \rtrw, \mathbf{F}\right) \mid s\in [0,\infty),\; \rtrw \in \R^d,\; \vc{F} \in \R^d \right\}
	\end{align*}
	is a measure--valued solution of the Aw-Rascle system \eqref{AR_multi} in $(0,T)\times \T^d$, with the bounded energy initial data \eqref{MV:IE:11} and dissipation defect $\mathcal{D}$,
	\begin{align*}
		\mathcal{D}\in L^{\infty}(0,T),\; \mathcal{D}\geq 0,
	\end{align*}
	if the following conditions hold:
	\begin{itemize}
		\item The unknowns enjoy the following regularity: 
		\begin{align}
			&  \vr=\lel \Nu_{t,x}; s  \ril  \in C_{\text{weak}}(0,T;L^{\gamma+1}(\T^d)),\label{reg1}\\
			& 	{\vr}^\alpha=\lel \Nu_{t,x}; {s}^\alpha  \ril  \in C_{\text{weak}}(0,T;L^{\frac{1}{\alpha}(\gamma+1)}(\T^d))  \text{ with } 0<\alpha< \gamma+1 ,\label{reg2}\\
			&\lel \Nu_{t,x};  \rtrw \ril \in L^\infty(0,T;L^2(\T^d)) \text{ with } \sqrt{\vr} \lel \Nu_{t,x};  \rtrw \ril =  \lel \Nu_{t,x}; \sqrt{s} \rtrw \ril,\label{reg3}\\
			&Q(\vr)=\lel \Nu_{t,x}; Q(s)  \ril \in L^2(0,T; W^{1,2}(\T^d)) \text{ with } \Grad Q(\vr)=\lel \Nu_{t,x}; \DQ  \ril.\label{reg4}
		\end{align}	
		\item For a.e. $\tau \in (0,T) $ and $\psi \in C^{1}([0,T]\times {\T^d})$ we have
		\begin{align} \label{mv eqn 1}
			\begin{split}
				&\int_{\T^d}\vr(\tau, \cdot) \psi(\tau, \cdot) \dx - \int_{ \T^d} \vr_0 \psi(0, \cdot) \dx \\
				&\quad = \int_{0}^{\tau} \int_{ \T^d} \left[ \vr  \partial_{t}\psi + \sqrt{\vr}\lel \Nu_{t,x};  \rtrw \ril \cdot \Grad \psi - \sqrt{\vr} \Grad Q(\vr) \cdot \Grad \psi \right] \dx \dt.
			\end{split}
		\end{align}
		\item There exists a defect measure $r^{M} \in L^\infty_{\text{weak-(*)}}(0,T;\mathcal{M}({\T^d};\R^{d\times d})) + \mathcal{M}([0,T]\times \T^d;\R^{d\times d})$ and $\xi\in L^{1}(0,T)$ such that for a.e. $\tau \in (0,T)$, for all $ \epsilon>0 $ and $\pmb{\phi}\in C^{1}([0,T]\times {\T^d};\R^d)$ we have	
		\begin{align}\label{mv eqn 2}
			\begin{split}
				&\int_{ \T^d} \sqrt{\vr}(\tau,x) \lel\Nu_{\tau,x}; \rtrw  \ril \cdot \pmb{\phi}(\tau,\cdot) \dx - \int_{ \T^d} \lel \Nu_{0} ;  \sqrt{s} \rtrw  \ril \cdot \pmb{\phi}(0,\cdot) \dx\\
				&= \int_{0}^{\tau} \int_{ \T^d} \sqrt{\vr}  \left[ \lel \Nu_{t,x};  \rtrw   \ril \cdot \partial_{t} \pmb{\phi} +  \lel \Nu_{t,x}; \rtrw \otimes  \rtrw  \ril : \Grad \pmb{\phi}- \lel \Nu_{t,x}; \rtrw \otimes \DQ \ril : \Grad \pmb{\phi} \right] \dx \dt\\
				&  + \lel r^M;\Grad \pmb{\phi}\ril_{\mathcal{M}([0,\tau]\times \T^d), C([0,\tau]\times \T^d)} 
			\end{split}
		\end{align}
	and
		\begin{equation}\label{mom-def}
		\left\vert \lel r^M;\Grad \pmb{\phi}\ril_{\mathcal{M}([0,\tau]\times \T^d), C([0,\tau]\times \T^d)} \right\vert \leq  \left( 1+\frac{1}{4 \epsilon} \right) \int_0^\tau \xi(t) \mathcal{D}(t) \Vert \pmb{\phi} \Vert_{C^1({\T^d})} \dt + \epsilon \Vert \pmb{\phi} \Vert_{C^1([0,\tau]\times {\T^d})} \mathcal{D}(\tau).
	\end{equation} 
		\item There exists a defect measure $ \mathcal{R} \in L^\infty_{\text{weak-(*)}}(0,T;\mathcal{M}(\T^d))+ \mathcal{M}([0,T]\times \T^d)$ such that for a.e. $ \tau \in (0,T) $ and $ \epsilon >0 $ we have
		\begin{align*}
			\begin{split}
				&\int_{ \T^d} \lel \Nu_{\tau ,x}; \frac{1}{2}  \levert \rtrw \rivert^2 + E(s) \ril \dx + \int_{0}^{\tau} 	\intO{ \lel \Nu_{t,x}; \vert \DQ \vert^2 \ril} \dt + \mathcal{D}(\tau) \\
				&\leq \int_{ \T^d} \lel \Nu_{0,x};  \frac{1}{2} \levert \rtrw \rivert^2 +E(s) \ril \dx + \int_{0}^{\tau} \intO{ \lel \Nu_{t,x}; \rtrw \cdot \DQ \ril }  \dt +\int_{(0,\tau)\times \T^d}{\text{d} \mathcal{R}}  ,
			\end{split}
		\end{align*}
		where
  \begin{align}\label{En:def}
      \left\vert \int_{(0,\tau)\times \T^d}{\text{d} \mathcal{R}}   \right\vert \leq C \left(1+\frac{1}{4 \epsilon} \right) \int_{0}^\tau \mathcal{D}(t) \dt + \epsilon \mathcal{D}(\tau).
  \end{align}
	Here we can identify $ \lel \Nu_{t,x}; E(s) \ril= E(\vr) $. 
	\end{itemize}
	\end{df}

	\begin{rmk}
	Our definition of measure-valued solution uses less regularity than a-priori available for $\sqrt{\vr}\vw$ in \eqref{est_mom}.
	    In fact, following our approximation procedure, we are only able to justify the uniform  $L^1((0,T)\times \T^d)$ bound on $\vr\vw\otimes\vw$. Therefore we also need to include possible concentration defect $r^M$ associated with this term in the definition of the solution \eqref{mv eqn 2}.
	\end{rmk}
	
	\begin{rmk}
    The regularity assumptions \eqref{reg1} and \eqref{reg2} imply
    $\lel \Nu_{t,x}; {s}^\alpha  \ril=\lel \Nu_{t,x}; {s}  \ril^\alpha$ for  a.e. $(t,x) \in (0,T) \times \T^d$. Recall that for a strictly convex function the Jensen inequality is an equality if and only if the measure is a Dirac delta.
    Choosing $\alpha=1+\frac{\gamma}{2}$ we get that the projection of Young measure on the first component is a Dirac delta. \\
    Therefore
    \[ \Nu_{t,x}= \delta_{\{\vr(t,x)\}}\otimes Y_{t,x} \quad \text{for  a.a.\ } (t,x) \in (0,T) \times \T^d, \] 
    where $Y \in L^{\infty}_{\text{weak-(*)}} \big( (0,T)\times \T^d ;\mathcal{P}( \R^d \times \R^{d})\big)$. 
\end{rmk}

\subsection{The main results}

The main results of the paper are formulated in the two theorems below.
\begin{thm}\label{Theorem:1}
Let $p(\vr)=\vr^\gamma $, $ \gamma \geq1$, and let the initial data satisfy \eqref{MV:IE:11}. Then, for any $ T>0 $,  there exists a measure-valued solution to \eqref{PM} in the sense of 
 Definition \ref{MVdef}.
\end{thm}

\begin{rmk}
In a more recent development, Abbatiello, Feireisl and Novotn\' y \cite{AFN21} and Wo\'znicki \cite{Wo} simplified the notion of measure-valued solutions to the so-called dissipative solutions (with no Young measure), for both compressible and incompressible Navier-Stokes systems with a general viscosity law describing non-Newtonian fluids. The same can also be used for the compressible Euler system, see Breit, Feireisl, and Hofmanov\'a \cite{BrFeHo}. This definition relies on writing the nonlinear term $ \vr \vu \otimes \vu $ in terms of variables $ (\vr, \vc{m}=\vr \vu) $, and on using the convexity of the function 
$$ (\varrho,\vc{m}) \mapsto \mathbf{1}_{\{\vr > 0\}}\frac{\vc{m} \otimes\vc{m} }{\vr} \colon (\xi \otimes \xi) =\mathbf{1}_{\{\vr > 0\}} \frac{\vert \vc{m} \cdot \xi \vert^2}{\vr}\quad \mbox{for all} \ \xi(\neq 0) \in \R^d .$$ 
However, in our case, we have two nonlinearities. The first one $ (\vr,\vw) \mapsto \vr \vw \otimes \vw $ is of course similar if we write it in terms of variable $ \sqrt{\vr} \vw $ or $ (\vr, \vc{m}=\vr \vw) $, but the nonlinearity $ (\Grad Q(\vr),\sqrt{\vr}{\vw}) \mapsto \sqrt{\vr} \vw \otimes \Grad Q(\vr) \colon (\xi \otimes \xi) $ is not convex. This is the main obstacle that prevents us from getting rid of the Young measure.
\end{rmk}

\begin{thm}\label{Theorem:2}
	Let  $ (\Nu, \mathcal{D}) $ be a measure-valued solution of \eqref{PM} in the sense of Definition \ref{MVdef}. Let $ (\bvr,\bvw) $ be a strong solution to the same system, s.t.
	\begin{equation*}
		\bvr \in C^1([0,T]; C^2(\T^d)),\; \bvw \in C^1([0,T]; C^2(\T^d; \R^d)) \quad
		\text{ with } \bvr >0,
		\end{equation*}
		 emanating from the initial data $ (\bvr_0,\bvw_0) \in (C^2(\T^d), C^2(\T^d;\R^d)) $ s.t. $ \bvr_0>0 $. If the initial data coincide, i.e.
\begin{equation*}
	\Nu_{0 ,x}= \delta_{\{ \bvr_0(x), \bvw_0(x)\}}, \text{ for a.e. }x\in \T^d
	\end{equation*} 
then $ \mathcal{D}=0 $, and 
\begin{align*}
		\Nu_{\tau ,x}= \delta_{\{ \bvr(\tau,x), \sqrt{\bvr}\bvw(\tau, x), \Grad Q(\bvr)(\tau,x) \}}, \text{ for a.e. }(\tau,x)\in(0,T)\times \T^d.
\end{align*}
\end{thm}
\begin{rmk}
   Theorem \ref{Theorem:2} asserts that, starting from the same initial data, the measure-valued solutions and the strong solution coincide in the  time interval $ [0,T)$ in which the strong solution exists. The existence of the strong solution is however not yet known, it is the subject of our future studies.
\end{rmk}	
 
\begin{rmk}
In the course of the proof it will become clear that for the Theorem \ref{Theorem:2} to hold,  the initial density needs to be strictly positive. In the context of weak-strong uniqueness results this assumption can sometimes be lifted, see for example \cite{FN2021}.
\end{rmk}	
The outline of the paper is as follows. Section \ref{Sec:3} is devoted to the formulation of the approximate problem. In Section \ref{existance}, we perform the limit passages leading ultimately to the measure-valued solution of the target system \eqref{PM} as stated in Theorem \ref{Theorem:1}. Section \ref{rel-ent-motivation} deals with the informal derivation of the relative energy inequality, justified only for smooth solutions and test functions.
And finally, in Section \ref{Sec:6}, we extend the relative energy functional to the measure-valued solutions and conclude the proof of Theorem \ref{Theorem:2}.

\section{Existence of solutions: the approximation scheme}\label{Sec:3}
At first, we fix positive parameters $\ep > 0$, $\delta > 0$ and $\kappa>0$.
\begin{itemize}
    \item For $\kappa >0$, we denote
\eqh{
\hat{\vr}_\kappa=\frac{\sqrt{\vr^2+\kappa^2}}{1+\kappa\sqrt{\vr^2+\kappa^2}} .}
\item The operator $f\to [f]^{\kappa}$, $\kappa=(\kappa_{t},\kappa_{x})$ is the standard smoothing operator, that applies  to the variables $x$ and $t$. \item We choose a system of functions $\{ \vc{w}_j \}_{j=1}^\infty \subset \DC(\Omega; R^d)$ that form an orthonormal basis of the space $L^2(\Omega; R^d)$. For $ n\in \mathbb{N}$, we define 
\[\ X_n = {\rm span} \left\{ \vc{w}_j \ \Big| \ 1 \leq j \leq n \right\} \]
with the scalar product
$$\langle\vw,\vc{v}\rangle=\intO{\vw\cdot\vc{v}},\quad \forall\vw,\vc{v}\in X_{n}.$$
\end{itemize}
Now, for fixed $ \kappa, n, \delta, \ep $, our basic level of approximation consists of:
\begin{itemize}
    \item {{the augmented quasi-linear continuity equation:}}
    \eq{
\pt \vr+\Div(\vr \vw)-\Div\lr{\lr{\ep+[\hat{\vr}_{\kappa} p'(\hat{\vr}_{\kappa})]^\kappa}\Grad\vr}=0,\label{AR1_app}}
\item {the Faedo-Galerkin approximation for the weak formulation of the momentum equation:}
\eq{\label{AR2_app}
\intO{\vr\vw(T)\cdot\pmb{\phi}}-\intO{\vc{m}_{0} \cdot\pmb{\phi}}+\delta\intTO{\Grad\vw :\Grad \pmb{\phi}}-\intTO{(\vr{\vw}\otimes{\vw}):\Grad \pmb{\phi}}\\+\intTO{\lr{\lr{\ep+[\hat{\vr}_{\kappa} p'(\hat{\vr}_{\kappa})]^\kappa}\Grad\vr\otimes\vw}:\Grad \pmb{\phi}}=0,
}
satisfied for any test function $\pmb{\phi}\in X_{n}$.
\end{itemize}
 In order to get our desired results we have to perform a limit passage in each of these parameters. This will be discussed in detail in Section \ref{existance}.
 Note that one of the construction levels involves viscous regularization of the equation for $\vw$. At this stage, the approximate system might be viewed as a pressureless version of Brenner's model introduced in \cite{B1,B2}. The existence result for this system with barotropic pressure was provided by Feireisl and Vasseur in \cite{FV}.\\
We consider the following approximation of the initial data for\eqref{AR1_app} and \eqref{AR2_app} 
\eq{\label{initial_delta}
&\vr_\delta(0,x)=\vr_{0,\delta}(x) \in C^{2,\nu}(\T^d),\quad (\vr\vw)_\delta(0,x)=\vc{m}_{0,\delta} (x)= \vc{m}_0(x),\\
&0<\underline{\vr_0} \leq\vr_{0,\delta}\leq\overline{\vr_0}<\infty,
}
such that  \eqref{MV:IE:2} is satisfied uniformly with respect to all approximation parameters,
and 
\[ \vr_{0,\delta} \rightarrow \vr_0 \in L^{\gamma+1}(\T^d) \text{ as } \delta \rightarrow 0. \]

\par

\begin{thm}\label{thap}
Let $\kappa,n,\ep,\delta$ be fixed positive parameters.
The approximate problem (\ref{AR1_app}-\ref{AR2_app}) admits a strong solution $\{\vr,\vw\}$ which enjoys the following regularity:
$$\vw\in C^{1}([0,T],X_{n}),$$
\eq{\label{class_vr}
\vr\in C([0,T];C^{2,\nu}(\T^d))\cap C^1([0,T]\times \T^d),\quad
\pt\vr\in C^{0,\nu/2}([0,T];C(\Omega)),}
\eqh{\underline{\vr}\leq\vr(t,x)\leq \overline{\vr},}
where ${0< \nu <1} $, and $\underline{\vr}, \overline{\vr}$ depend only on $\|\vw\|_{ C^{1}([0,T],X_{n})}$. 

\end{thm}

\pf The approximate system \eqref{AR1_app} and \eqref{AR2_app} is similar to the classical approximations used in the construction of solutions to compressible Navier-Stokes equations. Some of the basic arguments developed for this system are described in detail in the book of Novotn\'y and Straskraba \cite{ NS }, Feireisl and Novotn\'y \cite{ FN }, or in the pioneering work of Lions \cite{PLL}. The only differences in our approximation are 1) the additional quasi-linear diffusion in the continuity equation, and 2) the additional drift term in the momentum equation instead of the pressure. In the proof below we focus on explaining how to incorporate these differences in the classical proof.

\textbf{Step 1.} We set $\vw\in C([0, T ];X_{n})$ for which we find the mapping $\vw\rightarrow \vr(\vw)$   determining the unique solution to the continuity equation. Note, that our approximation of continuity equation includes artificial linear diffusion with parameter $\ep$ as well as truncation and regularisation of nonlinear diffusion coefficient. This resembles the approximation often used for the temperature equation in the full compressible Navier-Stokes(-Fourier) system, see Section 3.4.2 in \cite{FN}.  To prove the existence of strong solution to this equation we use the theory of  Ladyzenskaja, Solonnikov and Uralceva \cite{LSU}, developed in  the framework of continuously differentiable functions, and also for the Sobolev spaces. For convenience, we refer the reader to the book of Novotn\'y and Feireisl \cite{FN}, where the relevant theorem has been recalled as Theorem 10.24. Applying this theorem, we get that \eqref{AR1_app} admits a unique classical solution $\vr=\vr_\kappa$ in the class \eqref{class_vr}.

Moreover, using the structure of the equation \eqref{AR1_app} we can show the following lemma
\begin{lemma}\label{Lemma:approx}
Let $\vw\in C([0, T ];X_{n})$ and let $\vr(0,x)=\vr_{0,\delta}(x)$ be such that
\eqref{initial_delta} is satisfied.
Then problem \eqref{AR1_app}  admits  at most one strong solution $\vr$ in the class \eqref{class_vr}. This $\vr$ satisfies the following estimates:
\eq{\label{bounds}
&\|\vr^{-1}\|_{L^\infty((0,T)\times\T^d)}+\|\vr\|_{L^\infty((0,T)\times\T^d)}\leq C\\
&{\rm ess}\sup_{t\in(0,T)}\|\vr\|^2_{W^{1,2}(\T^d)}
+\intT{\lr{\|\pt\vr\|^2_{L^2(\T^d)} +\|\lap {\cal K}(\vr)\|^2_{L^2(\T^d)}}}\leq C,
}
where we denoted
\eqh{
{\cal K}_\kappa(\vr)=\int_{1}^\vr\lr{\ep+[\hat z_\kappa p'(\hat z_\kappa)]^\kappa}{\rm d}z.}
Moreover, the mapping $\vw\to\vr(\vw)$ maps bounded sets in $C([0, T ];X_{n})$ into the bounded sets of the space specified by the above regularity class, and is continuous with values in $C^1([0,T]\times(\T^d))$. 
\end{lemma}

\textbf{Step 2.} In this step, for sufficiently small time interval $[0,T]$, we find the unique solution to the momentum equation \eqref{AR2_app} applying the Banach fixed point theorem. Then we extend the existence result for the maximal time interval.\\
Let us prove that there exists $T=T(n)$ and $\vw\in C([0,T];X_{n})$ satisfying (\ref{AR2_app}).
To this purpose, we  apply the fixed point argument to the mapping
$${\cal{T}}:C([0,T];X_{n})\rightarrow C([0,T];X_{n}),$$
\begin{equation}\label{T}
{\cal{T}}[\vw](t)= {\cal{M}}\left[\vr(t),P_{n}\vc{m}_0+\int_{0}^{t}{P_{n}{\cal {N}}(\vw)(s) {\rm d}s}\right],
\end{equation}
where $P_{n}$ is the orthogonal projection of $L^{2}(\Omega)$ onto $X_{n}$,
$${\cal {N}} (\vw) =-\Div (\vr_\kappa \vw \otimes \vw) +\Div (\vr\Grad p(\vr)\otimes\vw)+\delta\lap\vw$$
and
$${\cal{M}}\left[\vr(t),\cdot\right]:X_{n}\rightarrow X_{n},\intO{\vr(t){\cal{M}}\left[\vr(t),\vc{w}\right]\phi}=\langle\vc{w},\pmb{\phi}\rangle, \quad\vc{w},\pmb{\phi}\in X_{n}.$$
First, let us observe that $P_n{\cal N}(\vw)(t)$ is bounded in $X_{n}$ for $t\in[0,T]$.
Using the equivalence of norms on the finite dimenssional space $X_n$ we can easily check that
\begin{align*}
    \|P_n{\cal N}(\vw)\|_{X_n}\leq C\left(\|\vw\|_{X_n}+\|\vr\|_{L^\infty(\T^d)}\|\vw\|_{X_n}^2+
\|\vr\|_{L^\infty(\T^d)}^{\gamma+1}\|\vw\|_{X_n}
\right).
\end{align*}
Concerning the operator ${\cal M}$, it is easy to see that provided $\vr(t,x)\geq\Un{\vr}>0$ we have
$$\|{\cal M}[\vr]\|_{{\cal L}(X_n,X_n)}\leq\Un{\vr}^{-1}.$$
Moreover, since ${\cal M}[\vr]-{\cal M}[\vr']={\cal M}[\vr']\left({\cal M}^{-1}[\vr']-{\cal M}^{-1}[\vr]\right){\cal M}[\vr]$
 for any $(t,x)\in([0,T]\times\Omega)$ we have that
$$\|{\cal M}[\vr]-{\cal M}[\vr']\|_{{\cal L}(C_n,X_n)}\leq c\Un{\vr}^{-2}\|(\vr-\vr')(t)\|_{L^1(\Omega)}$$
for $t\in[0,T]$. Thus by virtue of continuity of mapping $\vw\to\vr(\vw)$ and the estimates established in Lemma \ref{Lemma:approx},
one can verify that ${\cal{T}}[\vw]$ maps the ball
$$B_{R,T}=\left\{\vw\in C([0,T],X_{n}):\|\vw\|_{C([0,T],X_{n})}\leq R, \vw(0,x)=P_{n}\left(\frac{\vc{m}_{0}}{\vr_{0,\delta}}\right)\right\}$$
into itself, and it is a contraction for sufficiently small $T>0$. It therefore posses the unique fixed point satisfying \eqref{AR2_app} on the time interval $[0,T]$. In view of previous remarks, the proof of this step can be done by a minor modification of the procedure described in \cite{NS}, Section 7.7, so we skip this part.\\ 

\textbf{Step 3.} Additionally, the time regularity of $\vw$ may be improved directly by differentiating (\ref{T}) with respect to time and estimating the norm of the resulting right hand side in $X_{n}$, so we get
$$\vw\in C^{1}([0,T],X_{n}).$$

This is the crucial information that enables to extend our solution to the maximal time interval $[0,T]$, for $T$ arbitrary large, see again \cite{NS} for more details. $\Box$

\section{Passage to the limit with approximation parameters}\label{existance}
From this point on the original system \eqref{AR_multi} might be recovered through limit passages with the approximation parameters in the following order: $\kappa\to0$, $n\to\infty$, $\ep\to 0$, $\delta\to 0$.

The first three limit passages can be made within the standard weak formulation, and only the last limit passage will require switching to the definition of measure-valued solutions.

\subsection{Passage to the limit $\kappa\to 0$}
From now on, we denote the constructed solutions by $\{(\vr_\kappa, \vw_\kappa)\}_{\kappa > 0} $ to indicate that all parameters except for $ \kappa $ are fixed. Recall that the bounds \eqref{bounds} are uniform with respect to $\kappa$ provided that the dimension of Faedo-Galerkin approximation $n$ is fixed. They can be therefore used in order to deduce strong convergence of the density, we have
\eq{\label{cro_kappa}
\vr_\kappa\to\vr\quad \mbox{a.e.\ in} \ (0,T)\times\Omega.}
The classical energy estimate, on the other hand, yields the uniform estimates for the velocity vector field. To obtain 
it we first differentiate \eqref{AR2_app} with respect to time, and then observe that $\vw_{\kappa}$ can be used in it as a test function, we obtain
\eq{\label{energy_w_kappa}
\Dt\intO{\vr_\kappa|\vw_\kappa|^2}+\delta\intO{|\Grad\vw_\kappa|^2}=0.}
From this, integrating over time we obtain the following uniform estimates
\eq{\label{bounds_w}
&\|\vr_\kappa|\vw_\kappa|^2\|_{L^\infty(0,T;L^1(\T^d))}\leq C,\\
&\delta^\frac12\|\Grad\vw_\kappa\|_{L^2(0,T; L^2(\T^d))}\leq C,
}
where $ C $ is a generic constant, that is uniform with respect to $ \kappa $. Therefore, we obtain that, up to a subsequence, 
\eq{\label{cw_kappa}
\vw_\kappa\to\vw\quad\mbox{weakly in }\quad L^2(0,T; W^{1,2}(\T^d)).}
The uniform bounds \eqref{bounds} and \eqref{bounds_w}, the strong convergence of $\vr_\kappa$ and weak convergence of $\vw_\kappa$ allow us to pass to the limit in \eqref{AR1_app}. As a consequence, we have
\eq{
\pt \vr+\Div(\vr \vw)-\Div\lr{\lr{\ep+\vr p'(\vr)}\Grad\vr}=0,\quad \mbox{a.e. in }(0,T)\times\T^d.\label{AR1_kappa}}
We can also pass to the limit in all terms of \eqref{AR2_app} except for the convective term which is nonlinear in terms of $\vw$. To justify the limit passage in the convective term, let us first observe that the convergence of the momentum can be strengthened to
$$\vr_\kappa\vw_\kappa\to\vr\vw \quad\mbox{weakly}^{*}\mbox{ \ in} \ L^\infty(0,T;L^{2}(\T^d)),$$
due to the uniform estimates \eqref{bounds}, \eqref{bounds_w} and the strong convergence of the density. Next, one can show that for any $\pmb{\phi}\in \cup_{n=1}^\infty X_n$ the family of functions $\intO{\vr_\kappa\vw_\kappa(t)\cdot\pmb{\phi}}$ is bounded and equi-continuous in $C[0,T]$. Thus, via the Arzel\`a-Ascoli theorem and density of smooth functions in $L^2(\T^d)$ we get that
$$\vr_\kappa\vw_\kappa\to\vr\vw \quad\mbox{in} \ C_{weak}([0,T];L^{2}(\T^d)).$$
Finally, by the compact embedding $L^2(\T^d)\subset W^{-1,2}(\T^d)$ and the weak convergence of $\vw_\kappa$ \eqref{cw_kappa} we verify that
\eq{\label{conv_convective}
\vr_\kappa\vw_\kappa\otimes\vw_\kappa\to\vr\vw\otimes\vw\quad\mbox{weakly\ in} \ L^2((0,T)\times\T^d).}
With this convergence at hand, we can pass to the limit in \eqref{AR2_app} to obtain
\eq{\label{AR2_kappa}
\intO{\vr\vw(T)\cdot\pmb{\phi}}-\intO{\vc{m}_{0}\cdot\pmb{\phi}}+\delta\intTO{\Grad\vw :\Grad \pmb{\phi}}-\intTO{(\vr{\vw}\otimes{\vw}):\Grad \pmb{\phi}}\\+\intTO{\lr{\lr{\ep+\vr p'(\vr)}\Grad\vr\otimes\vw}:\Grad \pmb{\phi}}=0,
}
satisfied for any test function $\pmb{\phi}\in X_{n}$.

Let us also observe, that using lower semi-continuity argument, we can pass to the limit in the energy equality \eqref{energy_w_kappa} integrated w.r.t. time. It will now change into inequality, we obtain
\eq{\label{energy_w_n}
\intO{\vr|\vw|^2(T)}+\delta\intTO{|\Grad\vw|^2}\leq \intO{\vr_0|\vw_0|^2}.
}
Note that we could also write the above inequality for a.a. $\tau\in(0,T)$ instead of $T$.
\subsection{Passage to the limit $n\to \infty$} 
At this level we are not able to obtain a uniform upper and lower bound of density with respect to $n$ as in \eqref{bounds}. We have only $\vr_n\geq0$, and the uniform estimates of the velocity \eqref{bounds_w}. However, for $n$ fixed, we can still multiply  \eqref{AR1_kappa} by $p(\vr_n)$. Integrating by parts we obtain
\eqh{
\intO{\pt\vr_n p(\vr_n)}-\intO{\vr_n\vw_n\cdot\Grad p(\vr_n)}+\ep\intO{p'(\vr_n)|\Grad\vr_n|^2}+\intO{\vr_n|\Grad p(\vr_n)|^2}=0.
}
Seeing that $E'(\vr)=p(\vr)$,  we get
\eq{\label{energy_eq}
\Dt\intO{ E(\vr_n)}-\intO{\vr_n\vw_n\cdot\Grad p(\vr_n)}+\ep\intO{p'(\vr_n)|\Grad\vr_n|^2}+\intO{\vr_n|\Grad p(\vr_n)|^2}=0.
}
To estimate the second term on the l.h.s. we use some small constant $\sigma<<1$, and split
\eq{\label{split}
\left|\intO{\vr_n\vw_n\cdot\Grad p(\vr_n)}\right|\leq C_\sigma\intO{\vr_n|\vw_n|^2}+\sigma\intO{\vr_n|\Grad p(\vr_n)|^2}.
}
Integrating \eqref{energy_eq} w.r.t. time and using \eqref{bounds_w} to bound the first term on the r.h.s. of \eqref{split}, we finally obtain
\eqh{
\intO{ E(\vr_n)(T)}+\ep\intO{p'(\vr_n)|\Grad\vr_n|^2}+(1-\sigma)\intO{\vr_n|\Grad p(\vr_n)|^2}\leq C,
}
where the constant $C$ depends only on the initial data $ E_{0} $ and $\ep$.
Summarising the so-far obtained a-priori estimates, we get
\eq{\label{bound0}
&\Vert \sqvr_n \vw_n \Vert_{L^\infty(0,T;L^2( \T^d))} \leq C,\\
	& \delta^\frac12\Vert \Grad \vw_n \Vert_{L^2(0,T;L^2(\T^d))}\leq C,
	\\
	&\Vert E(\varrho_n) \Vert_{L^\infty(0,T;L^1( \T^d))} \leq C,\\
	& \ep^\frac12\Vert\sqrt{p'(\vr_n)}\Grad\vr_n \Vert_{L^2(0,T;L^2(\T^d))}\leq C,\\
	& \Vert\sqrt{\vr_n}\Grad p(\vr_n) \Vert_{L^2(0,T;L^2(\T^d))}\leq C.
}
Recalling the definition of $Q(\vr_n)$ \eqref{Q-defn}, we obtain in particular that
\begin{align}\label{bound1}
	\begin{split}
	& \Vert \Grad Q(\varrho_n) \Vert_{L^2(0,T;L^2(\T^d))}\leq C,\\
	&  \Vert Q(\varrho_n) \Vert_{L^\infty(0,T;L^{\sfrac{\lr{\gamma+1}}{\lr{\gamma+\sfrac{1}{2}}}}( \T^d))}\leq C,
\end{split}
\end{align}
which, due to the Poincar\'e inequality
\eqh{
\| Q\|_{L^2(\T^d)}^2
\leq C(d)\|\Grad Q\|_{L^2(\T^d)}^2
+\lr{\intO{Q}}^2
}
implies that 
\eqh{
 \Vert Q(\varrho_n) \Vert_{L^2(0,T;W^{1,2}(\T^d))}\leq C.
}
We now use the structure of $E$ and $Q$, along with the Sobolev imbedding  to write the above estimates in terms of $\vr_n$. In the more restrictive 3-dimensional case we obtain
\eq{\label{uniform_rho1}
& \Vert \varrho_n \Vert_{L^\infty(0,T;L^{\gamma+1}( \T^d))}\leq C,\\
&\Vert \varrho_n \Vert_{L^{(2\gamma+1)}(0,T; L^{6\gamma+3}(\T^d))}\leq C,
}
and so, by interpolation, we finally obtain 
\eq{\label{uniform_rho2}
\|\vr_n\|_{L^{4\gamma+3}(0,T; L^{4\gamma+3}(\T^d))}\leq C.
}
Note that for $d=2$ we get even higher integrability in space for $ \vr_n$. With this estimate at hand we again multiply the approximate continuity equation \eqref{AR1_kappa}, this time by $\vr_n$, we get
\eq{\label{rho2a}
\Dt\intO{\frac{\vr_n^2}{2}}
+\intO{\lr{\ep+\vr_n p'(\vr_n)}|\Grad\vr_n|^2}=
\intO{\vr_n\vw_n\cdot\Grad\vr_n}
.}
Notice that the time integral of the r.h.s. might be estimated as follows
\eq{\label{rho2b}
\left|
\intTO{\vr_n\vw_n\cdot\Grad\vr_n}\right|=\left|\frac12\intO{\vr_n^2 \, \Div\vw_n}
\right|\leq C\|\Div\vw_n\|_{L^2(0,T; L^2(\T^d))}\|\vr_n\|_{L^4(0,T; L^4(\T^d))}^2,
}
and the last term is bounded due to \eqref{uniform_rho2}. Thanks to this estimate, we obtain additionally that
\begin{align}\label{bound2}
	\begin{split}
		& \ep^\frac12\Vert \Grad \vr_n \Vert_{L^2(0,T;L^2(\T^d))}\leq C,\\
		& \ep^\frac12\Vert\sqrt{p'(\vr_n)}\Grad\vr_n \Vert_{L^2(0,T;L^2(\T^d))}\leq C.
\end{split}
\end{align}
In order to pass to the limit $n\to+\infty$, we must first show that the sequence $\{\vr_n\}_{n\in \mathbb{N}}$ converges strongly to $\vr$. 
To this purpose we use the following generalization of the classical Aubin-Lions lemma.
\begin{lemma}\label{Lemma:strong}
Let
\eq{\label{time_der}
\sup_{n\in\mathbb{N}}\|\pt\vr_n\|_{L^2(0,T; W^{-m,2}(\T^d))}\leq C_1<+\infty
}
for some $m\in \mathbb{N}$, and let
\eq{\label{ass:Q}
\sup_{n\in\mathbb{N}}\| Q(\vr_n)\|_{L^2(0,T; W^{1,2}(\T^d))}\leq C_2<+\infty.
}
Then, up to extraction of a subsequence, we have
\eqh{
\vr_n\to\vr\quad \mbox{in}\quad {L^2((0,T)\times\T^d)}.
}
\end{lemma} 
\pf {\bf Step 1.} First note that the assumption \eqref{time_der} is satisfied uniformly w.r.t. $n$. For this, we estimate the time derivative from \eqref{AR1_kappa} using \eqref{bound0}, \eqref{uniform_rho1} and \eqref{bound2}.
Indeed, we have
\begin{align*}
    \Vert \vr_n \vw_n \Vert_{L^\infty(0,T;L^{\sfrac{2(\gamma+1)}{(\gamma+2)}}( \T^d))}+\Vert \vr_n \Grad p(\vr_n) \Vert_{L^2(0,T;L^{\sfrac{2(\gamma+1)}{(\gamma+2)}}( \T^d))} \leq C.
\end{align*}
We know that the spaces $W^{-m,2}(\T^d),\ L^2(\T^d), $ and $W^{1,2}(\T^d)$ have common basis $\{v_i\}_{i=1}^{+\infty}$ (orthonormal in $L^2(\T^d)$ and orthogonal in $W^{1,2}(\T^d)$ and in $W^{-m,2}(\T^d)$). 
Note that in such a case also $\{\Grad v_i\}_{i=1}^{+\infty}$ is an orthogonal system in
$L^2(\T^d)$ (but not a basis).
Such functions can be constructed as tensor products of Fourier's one dimensional basic functions i.e. constant function, $\{{\rm sin}(2j\Pi x)\}_{j=1}^\infty$ and  $\{{\rm cos}(2j\Pi x)\}_{j=1}^\infty$.\\
Let us denote $\lambda_i:=\frac{\|\Grad v_i\|_{L^2(\T^d)}}{\|v_i\|_{L^2(\T^d)}}$.
The compactness of the embedding $W^{1,2}\subset L^2$ implies that $\lambda_i \rightarrow \infty$ as $i \rightarrow \infty$. Changing the order if necessary we can regard $\{\lambda_i\}_{i=1}^\infty$ as monotone non-decreasing sequence.

Let $\{w_j\}_{j=1}^\infty$ be a orthonormal basis of $L^2(0,T)$.
We denote ${V}_k={\rm{lin}}\{\{w_j\cdot v_i\}_{j=1}^\infty\}_{i=1}^k$, 
and by $\Pi_{V_k}$ we denote the projection operator on this space.

The main idea of the proof is to replace $\vr_n$ by $\Pi_{V_k}\vr_n$, to which we will apply the classical Arzel\`a-Ascoli argument. 

{\bf Step 2.} Denote $\Ov\vr^\sigma:=\max\{\sigma,\vr\}$ for some small $\sigma>0$, we use it to split
\eq{\label{triangle}
&\|\vr_n-\Pi_{V_{k}}\vr_n\|_{L^2((0,T)\times\T^d)}\\
&\leq 
\|\vr_n-\Ov{\vr_n}^\sigma\|_{L^2((0,T)\times\T^d)}
+\|\Ov{\vr_n}^\sigma-\Pi_{V_{k}}\Ov{\vr_n}^\sigma\|_{L^2((0,T)\times\T^d)}
+\|\Pi_{V_{k}}\Ov{\vr_n}^\sigma -\Pi_{V_{k}}\vr_n\|_{L^2((0,T)\times\T^d)}.
}
Note that
\eqh{
\|\vr_n-\Ov{\vr_n}^\sigma\|_{L^2((0,T)\times\T^d)}\leq C\sigma
}
and
\eqh{
\|\Pi_{V_{k}}\Ov{\vr_n}^\sigma -\Pi_{V_{k}}\vr_n\|_{L^2((0,T)\times\T^d)}\leq C\sigma,}
where $C$ depends only on the measure of $\T^d$.
To treat the middle term of \eqref{triangle} we use that  \eqref{ass:Q} implies that
\begin{align*}
    \sup_{n\in\mathbb{N}}\levertl Q(\Ov{\vr_n}^\sigma)\rivertl_{L^2(0,T; W^{1,2}(\T^d))}\leq C,
\end{align*}
therefore
\eqh{
\sup_{n\in\mathbb{N}}\| \Grad \Ov{\vr_n}^\sigma\|_{L^2((0,T)\times\T^d)}=
\sup_{n\in\mathbb{N}}\levertl \frac{1}{Q'(\Ov{\vr_n}^\sigma)}\Grad Q(\Ov{\vr_n}^\sigma)\rivertl_{L^2((0,T)\times\T^d)}\leq \frac{C}{Q'(\sigma)}.
}
And so, we have
\eqh{
\|\Ov{\vr_n}^\sigma-\Pi_{V_{k}}\Ov{\vr_n}^\sigma\|_{L^2((0,T)\times\T^d)}\leq \frac{C}{\lambda_k\, Q'(\sigma)}.
}
Therefore, plugging it into \eqref{triangle}, we obtain
\eqh{\|\vr_n-\Pi_{V_{k}}\vr_n\|_{L^2((0,T)\times\T^d)}\leq C(k,\sigma),}
where $C=C(k,\sigma)\to 0$ as $k\to+\infty$ for $\sigma$ fixed.\\

{\bf Step 3.} Clearly, we have
\eqh{
\|\Pi_{V_k}\pt\vr_n\|_{L^2(0,T; W^{-m,2}(\T^d))}\leq \|\pt\vr_n\|_{L^2(0,T; W^{-m,2}(\T^d))}\leq C_1.
}
On the other hand, on the finite-dimensional space $V_k$ all the norms are equivalent. In particular, $\|\cdot\|_{W^{-m,2}(\T^d)}$ and $\|\cdot\|_{L^2(\T^d)}$ are equivalent, therefore
\eqh{
\|\pt\Pi_{V_k}\vr_n\|_{L^2((0,T)\times\T^d)}=\|\Pi_{V_k}\pt\vr_n\|_{L^2((0,T)\times\T^d)}
\leq \|\Pi_{V_k}\pt\vr_n\|_{L^2(0,T; W^{-m,2}(\T^d))}\leq C_1.
}
We can now apply the Arzel\`a-Ascoli theorem to prove the compactness of the sequence $\{\Pi_{V_k}\vr_n\}_{n\in \mathbb{N}}$ in $L^2((0,T)\times\T^d)$ as $k\to\infty$.
We conclude the proof by taking first $k\to\infty$ and then $\ep\to 0$ in Steps 2 and 3. $\Box$\\
With this lemma at hand, and with the uniform bounds \eqref{uniform_rho1} we deduce that, up to a subsequence
\eq{\label{strong_vr_n}
\vr_n\to\vr\quad\mbox{strongly in}\quad L^{p_1}(0,T; L^{q_1}(\T^d))\cap L^{p_2}(0,T; L^{q_2}(\T^d)),}
for all $1\leq p_1<\infty,\; 1\leq q_1 <\gamma+1,\;  1\leq p_2 < 2\gamma+1\leq  $ and $ 1\leq q_2< 6\gamma +3 $.
Having this, passage to the limit in the continuity equation leads to a weak formulation
\eq{\label{AR1_n}
 			&\int_{\T^d} \vr\psi(\tau, \cdot)  \psi(\tau, \cdot) \dx - \int_{ \T^d} \vr_0(\cdot)  \psi(0, \cdot) \dx \\
 			&\quad = \int_{0}^{\tau} \int_{ \T^d} \left[  \vr  \partial_{t}\psi + \vr\vw \cdot \Grad \psi - 
 			(\ep+\vr p'(\vr))\Grad\vr
 			\cdot \Grad \psi \right] \dx \dt,
}
satisfied for  a.e. $\tau \in (0,T) $ and every $\psi \in C^{1}([0,T]\times {\T^d})$.

\bigskip

The biggest problem in passing to the limit $n\to+\infty$ in the momentum equation is identification of the limit in the drift term $\vr_n\Grad p(\vr_n)\otimes\vw_n=\frac1\gamma\Grad E(\vr_n)\otimes\vw_n$. This can be sorted out thanks to the uniform estimate for $\Grad\vw_n$, the estimate for $\pt E(\vr_n)$ in some negative Sobolev space, and the compensated compactness type of argument. Let us recall the following lemma, whose proof can be found in \cite[Lemma 5.1]{PLL}.
\begin{lemma}\label{lem_compac} Let $g_n$, $h_n$ converge weakly to $g$, $h$ respectively in $L^{p_1}(0,T;L^{p_2}(\Omega))$, $L^{q_1}(0,T;L^{q_2}(\Omega))$ where $1\leq p_1,p_2 \leq \infty$, $\frac{1}{p_1}+ \frac{1}{q_1}=\frac{1}{p_2}+ \frac{1}{q_2}= 1$. We assume in addition that
\begin{itemize}
  \item $\partial_t g_n$ is bounded in $L^1(0,T;W^{-m,1}(\Omega))$ for some $m\geq 0$ independent of $n$.
  \item $\|h_n-h_n(t,\cdot+\xi)\|_{L^{q_1}(0,T;L^{q_2}(\Omega))} \rightarrow 0$ as $|\xi|\rightarrow 0$, uniformly in $n$.
\end{itemize}
Then $g_nh_n$ converges to $gh$ in $\mathcal{D}'$.
\end{lemma}
First note that $g_n=\Grad E(\vr_n)=\vr_n\Grad p(\vr_n)$ is uniformly bounded in $L^2(0,T; L^{\frac{2(\gamma+1)}{\gamma+2}}(\T^d))$, and $h_n=\vw_n$ is uniformly bounded in $L^2(0,T; L^6(\T^d))$, so the product of $g_n \cdot h_n$ is bounded in $L^1(0,T; L^1(\T^d))$ provided that $\gamma\geq \frac12$. Moreover, due to strong convergence of the density \eqref{strong_vr_n} and uniform bound for $\Grad E(\vr_n)$, we know that
\eqh{
\Grad E(\vr_n)\to \Grad E(\vr) \quad \mbox{weakly\ in}\quad L^2(0,T; L^{\frac{2(\gamma+1)}{\gamma+2}}(\T^d)). }
Next, using \eqref{def:E} and \eqref{AR1_kappa} we compute
\eq{\label{pt_E}
\pt E(\vr_n)=p(\vr_n)\pt\vr_n&=-\Div(\vr_n\vw_n)p(\vr_n)+\Div(\vr_n\Grad p(\vr_n))p(\vr_n)+\ep\Div(\Grad\vr_n)p(\vr_n)\\
&=-\Div(p(\vr_n)\vr_n\vw_n)+\vr_n\vw_n\cdot \Grad p(\vr_n)
+\Div(\vr_n p(\vr_n)\Grad p(\vr_n))-\vr_n|\Grad p(\vr_n)|^2\\
&\quad +\ep\Div(p(\vr_n)\Grad\vr_n)-\ep p'(\vr_n)|\Grad \vr_n|^2.
}
From the uniform estimates \eqref{bound0}, we check that the second, the fourth, and the sixth terms are bounded in $L^1(0,T; L^1(\T^d))$. The first, the third and the fifth terms are bounded in $L^1(0,T;W^{-1,1}(\T^d))$ provided that $\sqrt{\vr_n}p(\vr_n)\in L^1(0,T; L^2(\T^d))$ uniformly in $n$, which is true due to uniform bounds for $\vr_n$ \eqref{uniform_rho1}.

From what has been written we deduce that $\pt g_n$ is uniformly bounded in $L^1(0,T; W^{-2,1}(\T^d))$. We also know that the second condition in Lemma \ref{lem_compac} is satisfied because $\Grad\vw_n$ is uniformly bounded in $L^2(0,T; L^2(\T^d))$. Thanks to this, we can now verify that

\eq{\label{conv_product}
g_n\cdot h_n=\Grad E(\vr_n)\cdot\vw_n\to g\cdot h= \Grad E(\vr)\cdot\vw
}
in the sense of distributions. Moreover, we note that $ \Grad E(\vr)\cdot\vw = \sqrt{\vr}\Grad Q(\vr) \cdot \vw $.  The limit in the convective term can be proven exactly the same way as in \eqref{conv_convective}.

Letting $n\to+\infty$ in the approximate momentum equation \eqref{AR2_kappa}, we obtain:
\eq{\label{AR2_n}
&\intO{\vr\vw(\tau,\cdot)\pmb{\phi}(\tau,\cdot)}-\intO{\vc{m}_{0}(\cdot)\pmb{\phi}(0,\cdot)}\\
&=\int_0^\tau\intO{\left[\vr\vw\cdot \pmb{\phi}-\delta\Grad\vw :\Grad \pmb{\phi}+(\vr{\vw}\otimes{\vw}):\Grad \pmb{\phi}-\lr{\lr{\ep+\vr p'(\vr)}\Grad\vr\otimes\vw}:\Grad \pmb{\phi}\right]}\,\dt,
}
satisfied for  a.e. $\tau \in (0,T) $ and every $\pmb{\phi}\in C^{1}([0,T]\times {\T^d})$.
Moreover, we can also pass to the limit in the approximate energy estimate \eqref{energy_w_n} and integrated in time equation \eqref{energy_eq} that will now change into the inequality
\eqh{
\intO{ E(\vr)(T)}+\ep\intTO{p'(\vr)|\Grad\vr|^2}+\intTO{\vr|\Grad p(\vr)|^2}\\
\leq\intO{ E(\vr_0)}+\intTO{\vr\vw\cdot\Grad p(\vr)}.
}
The main difficulty in this step comes from the fact that the bounds \eqref{bound2} do no longer provide us with estimates independent of $\ep$. However, the arguments for the strong convergence of the density, as well as justification of convergence in the most restrictive nonlinear terms \eqref{conv_convective}, \eqref{conv_product} are the same as in the previous steps. Note, in particular, that the two last terms in \eqref{pt_E} can be estimated as previously. Indeed, we have
\eqh{
|\ep p(\vr_\ep)\Grad\vr_\ep|=\ep\vr_\ep^\gamma|\Grad\vr_\ep|=
\sqrt{\ep}\sqrt{\vr^{\gamma-1}}|\Grad\vr_\ep|\sqrt{\ep}\sqrt{\vr_\ep^{\gamma-1}}\vr_\ep
\leq C\ep_\ep p'(\vr_\ep)|\Grad\vr_\ep|^2+ C\ep p'(\vr_\ep)\vr_\ep^2.
}
Both of the terms are bounded at least in $L^1(0,T; L^1(\T^d))$  due to the uniform bounds \eqref{bound0}, \eqref{uniform_rho1}. 
Thanks to the strong convergence of the density and weak convergence of the velocity, passages to the limits in the convective terms and drift terms, can be done as on the previous level. Indeed for the latter we choose $g_\ep= \sqrt{\vr_\ep} \Grad Q(\vr_\ep) $ and $h_\ep= \vw_\ep$ and apply Lemma \ref{lem_compac}.

Next, one also needs to verify, that all $\ep$-dependent terms in the continuity equation, as well as the momentum equation converge to $0$ for $\ep\to0$. For this, we simply use the uniform estimate \eqref{bound2} thanks to which
\eqh{
\ep\Grad\vr_\ep\to 0\quad \mbox{strongly\ in\ }\ L^2(0,T;L^2(\T^d)).}
So, for example, the $\ep$-dependent term in the momentum equation \eqref{AR2_n} vanishes because
\eqh{
\left|\int_0^\tau\intO{\ep\Grad\vr_\ep\otimes\vw_\ep}\, \dt\right|
\leq \ep \|\vr_\ep\|_{L^2(0,T;L^2(\T^d))}\|\vw_\ep\|_{L^2(0,T;L^2(\T^d))}\to 0,}
when $\ep\to0$ thanks to the uniform bound for $\vw_\ep$ from \eqref{bound0}.\\
The summary of this section can be formulated as a separate result
\begin{lemma}\label{Lemma:delta}
Let $\delta>0$ be fixed, and let the initial data satisfy \eqref{initial_delta}.
Then, there exists a sequence of weak solutions $\{\vr_\delta,\vw_\delta\}_{\delta>0}$ such that:\\
1. The continuity equation
\eq{\label{AR1_delta}
 			&\int_{\T^d} \vr_\delta\psi(\tau, \cdot)  \psi(\tau, \cdot) \dx - \int_{ \T^d} \vr_{\delta,0}(\cdot)  \psi(0, \cdot) \dx \\
 			&\quad = \int_{0}^{\tau} \int_{ \T^d} \left[  \vr_\delta  \partial_{t}\psi + \vrd\vwd \cdot \Grad \psi - 
 			\vrd p'(\vrd)\Grad\vrd
 			\cdot \Grad \psi \right] \dx \dt
}
 is satisfied for  a.e. $\tau \in (0,T) $ and every $\psi \in C^{1}([0,T]\times {\T^d})$;\\
2. The momentum equation
\eq{\label{AR2_delta}
&\intO{\vrd\vwd(\tau,\cdot)\pmb{\phi}(\tau,\cdot)}-\intO{\vc{m}_{0}(\cdot)\pmb{\phi}(0,\cdot)}\\
&=\int_0^\tau\intO{\left[\vrd\vwd\cdot \pmb{\phi}-\delta\Grad\vwd :\Grad \pmb{\phi}+(\vrd\vwd\otimes\vwd):\Grad \pmb{\phi}-\sqrt{\vrd} \Grad Q(\vrd) \otimes \vwd : \Grad \pmb{\phi} \right]}\,\dt
}
is satisfied for  a.e. $\tau \in (0,T) $ and every $\pmb{\phi}\in C^{1}([0,T]\times {\T^d})$;\\
3. The energy estimates
\eq{\label{energy_w_delta}
\intO{\vr|\vwd|^2(\tau)}+\delta\int_0^\tau\intO{|\Grad\vwd|^2}\,\dt\leq \intO{\vr_{\delta,0}|\vw_{\delta,0}|^2},
}
and
\eqh{
\intO{ E(\vrd)(\tau)}+\int_0^\tau\intO{|\Grad Q(\vrd)|^2}\,\dt\leq\intO{ E(\vr_{\delta,0})}+\int_0^\tau\intO{\vwd\cdot \sqrt{\vrd} \Grad Q(\vrd)}\,\dt
}
are satisfied for  a.e. $\tau \in (0,T) $;
\end{lemma}

%
\subsection{Passage to limit $\delta\to 0$: Proof of Theorem \ref{Theorem:1}}
This section is dedicated to completion of the proof of Theorem \ref{Theorem:1}. 
First note that the sequence $ \{(\vrd, \vwd) \}_{\delta>0}$ constructed in Lemma \ref{Lemma:delta} satisfies the bounds
\eqref{bound0}, \eqref{bound1}, \eqref{uniform_rho1} uniformly in $\delta$.  
Note also that at this level we do not have any information on the bounds of $ \vwd $ independently of $\vrd$. However, we still can prove a strong convergence of the sequence approximating the density.
From the assumption  \eqref{MV:IE:2}, or   \eqref{MV:IE:11} equivalently, the initial data is uniformly bounded, and thus, uniformly w.r.t. $\delta$ we have
\begin{align*}
	\begin{split}
		&\Vert E(\varrho_\delta) \Vert_{L^\infty(0,T;L^1( \T^d))} \leq C,\\
		& \Vert \Grad Q(\varrho_\delta) \Vert_{L^2(0,T;L^2(\T^d))}\leq C,\\
		&\Vert \sqvr_\delta \vw_\delta \Vert_{L^\infty(0,T;L^2( \T^d))} \leq C,\\
		& \delta^{\frac{1}{2}} \Vert \Grad \vw_\delta \Vert_{L^2(0,T;W^{1,2}(\T^d))}\leq C,
	\end{split}
\end{align*}
where $C$ depends only on $ E_0 $. As a consequence we have 
\begin{align*}
	\begin{split}
		& \Vert \varrho_\delta \Vert_{L^\infty(0,T;L^{\gamma+1}(\T^d))}\leq C,\\
		&  \Vert Q(\varrho_\delta) \Vert_{L^\infty(0,T;L^{\sfrac{\lr{\gamma+1}}{\lr{\gamma+\sfrac{1}{2}}}}( \T^d))}\leq C,\\
		& \Vert  Q(\varrho_\delta) \Vert_{L^2(0,T;W^{1,2}(\T^d))}\leq C.
	\end{split}
\end{align*}
Thus we can apply Lemma \eqref{lem_compac} to deduce 
\begin{align*}
	&\varrho_\delta \rightarrow \varrho \text{ in } L^{p}(0,T; L^{q}(\T^d)),\quad \forall p<\infty,\ q<\gamma+1.
\end{align*}
From this we deduce
\begin{align*}
	& Q(\varrho_\delta) \rightarrow Q(\varrho) \text{ in } L^p(0,T; L^{q}(\T^d)),\; \forall p<\infty,\ q<\frac{\gamma+1}{\gamma+\sfrac12},\\ 
	& \Grad Q(\vr_\delta) \rightarrow \Grad Q(\vr) \text{ weakly in } L^2(0,T; L^{2}(\T^d)).
\end{align*}
On the other hand, we only know that
\begin{align*}
	&\sqrt{\vr_\delta} \vw_\delta \rightarrow \overline{\sqrt{\varrho} \vw  }\ \text{ weak-(*)ly in }  L^\infty(0,T; L^{2}(\T^d)).
\end{align*}
Using  the fundamental theorem of the Young measures (see M\'alek et al. \cite[Section 4.2, Theorem 2.1]{MNRR}), we can say that there exists a Young measure $ \Nu $ generated by $ \left\{(\vr_\delta, \sqrt{\vr_\delta} \vw_\delta,  \Grad Q(\vr_\delta))\right\}_{\delta>0} $. Due to Lemma \ref{Lemma:strong} and \cite[Corollary 3.4]{MPI} we get that the projection of the Young measure on the first component is a Dirac delta i.e., our measure is of a particular form:
\[ \Nu_{t,x}= \delta_{\{\vr(t,x)\}}\otimes Y_{t,x} \quad \text{for  a.a.\ } (t,x) \in (0,T) \times \T^d. \] 
Therefore, a correct phase space for the Young measure is
\begin{align*}
	\mathcal{F}= \left\{\left(s, \rtrw, \mathbf{F}\right) \mid s \in [0,\infty),\; \rtrw \in \R^d,\; \mathbf{F} \in \R^d \right\}.
\end{align*}
Clearly, we notice that 
\begin{align*}
	\begin{split}
	&\overline{\sqrt{\varrho} \vw  } = \lel \Nu_{t,x}; \rtrw \ril,\quad\vr=\lel \Nu_{t,x}; s  \ril, \\
	& Q(\vr)=\lel \Nu_{t,x}; Q(s)  \ril,\quad \Grad Q(\vr)=\lel \Nu_{t,x}; \DQ  \ril.
\end{split}
\end{align*}
This allows us to conclude that 
\begin{align*}
	& \varrho_\delta \vw_\delta \rightarrow \sqrt{\vr}\lel \Nu_{t,x};  \rtrw \ril\  \text{ weak-(*)ly\ in\ } L^\infty(0,T;L^{\frac{2(\gamma+1)}{(\gamma+2)}} (\T^d)),\\
	& \varrho_\delta \Grad p(\vr_\delta)  \rightarrow \sqrt{\vr} \Grad Q(\vr) \ \text{ weak-(*)ly\ in\ } L^\infty(0,T;L^{\frac{2(\gamma+1)}{(\gamma+2)}} (\T^d)).
\end{align*}
Thus we can  pass to limit in the continuity equation to obtain that for a.e. $\tau \in (0,T) $ and $\psi \in C^{1}([0,T]\times {\T^d})$, we have
\begin{align} \nonumber
	\begin{split}
		&\int_{\T^d}\vr(\tau, \cdot) \psi(\tau, \cdot) \dx - \int_{ \T^d} \vr_0 \psi(0, \cdot) \dx \\
		&\quad = \int_{0}^{\tau} \int_{ \T^d} \left[ \vr  \partial_{t}\psi + \sqrt{\vr}\lel \Nu_{t,x};  \rtrw \ril \cdot \Grad \psi - \sqrt{\vr} \Grad Q(\vr) \cdot \Grad \psi \right] \dx\, \dt.
	\end{split}
\end{align}

Our goal now is to identify the limits of the nonlinear terms in the momentum equation and energy inequality.  
At this point we give a lemma that generalizes the result by Feireisl et al., see \cite{FGSW}[Lemma 2.1].
\begin{lemma}\label{me-lemma}
Let $M,N \in \mathbb{N}$, $Q \subset \R^M$ be a bounded domain and $ \left\{ \vc{Z}_n \right\}_{n\in \mathbb{N}} $, $ \vc{Z}_n: Q\rightarrow \R^N $ be a sequence of equi-integrable functions which generates the Young measure $ \{\nu_y\}_{y\in Q} $. Let $ f,g: \R^N\rightarrow [0,\infty) $ be continuous functions such that 
\[ \sup_n (\Vert f(\vc{Z}_n ) \Vert_{L^1(Q}) + \Vert g(\vc{Z}_n ) \Vert_{L^1(Q) } < \infty \]
	and let $ F $ be a continuous function such that 
	\[ F:\R^N \rightarrow R,\; \vert F(\vc{z}) \vert \leq \epsilon f (\vc{z}) + \frac{1}{4 \epsilon} g(\vc{z}), \text{ for all } \vc{z} \in \R^N  \text{ and  for some }\epsilon>0. \]
	Denote the concentration defects as
	\[ F_\infty = \overline{F} - \lel \nu_y; F({\vc{Y}}) \ril ,\; f_\infty = \overline{f} - \lel \nu_y; f({\vc{Y}}) \ril  \text{ and }g_\infty = \overline{g} - \lel \nu_y; g({\vc{Y}}) \ril,  \]
	where $  \overline{F},  \overline{f}    $ and $  \overline{g}  $ are the weak-(*) limits of $ \{F( \vc{Z}_n)\}_{n\in \mathbb{N}} $,  $ \{f( \vc{Z}_n)\}_{n\in \mathbb{N}} $ and  $ \{g( \vc{Z}_n)\}_{n\in \mathbb{N}} $ in $ \mathcal{M}(\bar{Q}) $, respectively, and $\vc{Y}$ stands for a dummy variable. Then we have 
	\[ \vert F_\infty \vert \leq \epsilon f_\infty + \frac{1}{4\epsilon} g_\infty.\]
	\end{lemma}

From the uniform bounds we have
\begin{align*}
	&\frac{1}{2} \varrho_\delta \vert \vw_\delta \vert^2 + E(\vr_\delta) \rightarrow \bar{E} \ \text{ weak-(*)ly in } L^\infty_{\text{weak-(*)}}(0,T;\mathcal{M}(\T^d)),\\
	& \vert \Grad Q(\varrho_\delta) \vert^2 \rightarrow \sigma\  \text{ weak-(*)ly in } \mathcal{M}((0,T)\times \T^d).
\end{align*}
Applying Lemma \ref{me-lemma} we can show that
\begin{align*}
	&E_\infty = \bar{E} - \lel \Nu_{\tau ,x}; \frac{1}{2} \levert \rtrw \rivert^2+ E(s) \ril \in  L^\infty_{\text{weak-(*)}}(0,T;\mathcal{M}^{+}(\T^d)),\\
	&	\sigma_\infty = \sigma - \lel \Nu_{t,x}; \vert \DQ \vert^2 \ril\in  \mathcal{M}^{+}((0,T)\times \T^d).
\end{align*}
 We also define the total defect as
\begin{align}\label{defectDpf}
	\mathcal{D}(\tau) = E_\infty (\tau )[\T^d]+ \sigma_\infty[[0,\tau] \times \T^d].
\end{align}

We denote  the weak-(*) limit of the convective term $ \vr_\delta \vw_\delta \otimes \vw_\delta $ by $ \overline{\varrho \vw \otimes \vw} $. 
The defect measure associated to the convective term is
\[ r^{M,1} = \overline{\varrho \vw \otimes \vw} - \lel \Nu_{t,x}; \sqrt{s} \vc{v} \otimes \sqrt{s} \vc{v}\ril .\] 
We note that 
\[ \vert (\varrho_\delta \vw_\delta \otimes \vw_\delta)_{ij} \vert = \vert \vr_\delta w_{\delta_{i}} w_{\delta_{j}} \vert \leq \varrho_\delta  \vert \vw_\delta  \vert^2,\; \text{ for all } i,j=1,\cdots,d .\]
Therefore, using the Lemma \ref{me-lemma} in each entry of the matrix, we deduce that 
\[ \vert r^{M,1}_{ij} (\tau ) \vert \leq C E_\infty (\tau ),\; \text{ for all } i,j=1,\cdots,d .\]
We also have $ r^{M,1}_{ij}\in L^\infty_{\text{weak-(*)}}(0,T; \mathcal{M}(\T^d)) $.
Thus from \eqref{defectDpf}, we obtain 
\begin{align*}
	\left\vert \lel r^{M,1};\Grad \pmb{\phi}\ril_{\mathcal{M}([0,\tau]\times \T^d), C([0,\tau]\times \T^d)} \right\vert \leq  C  \int_0^\tau \Vert \pmb{\phi} (t) \Vert_{C^1(\T^d)} \mathcal{D}(t) \dt  .
\end{align*}
For the term $ \sqrt{\vr_\delta} \vw_\delta \otimes \Grad Q(\varrho_\delta) $, we consider the defect measure 
\[ r^{M,2} = \overline{\varrho \vw \otimes \Grad p(\vr)} - \lel \Nu_{t,x}; \rtrw \otimes \DQ \ril . \]
In this case, we observe that 
\[ \left\vert (\sqrt{\vr_\delta} \vw_\delta \otimes \Grad Q(\varrho_\delta) )_{ij} \right\vert = \vert \sqrt{\vr_\delta} w_{\delta{_i}} (Q(\varrho_\delta)_{x_j}) \vert \leq \frac{1}{4 \epsilon}\varrho_\delta \vert \vw_\delta \vert^2 + \epsilon \vert \Grad Q(\varrho_\delta ) \vert^2,\]
for $ \epsilon>0 $.
Applying Lemma \ref{me-lemma} again, we have 
\begin{align*}
	\levert r^{M,2}_{ij} \rivert \leq \frac{1}{4 \epsilon}  E_\infty  + \epsilon \sigma_\infty, 
\end{align*}
in $ \mathcal{M}([0,T]\times \T^d)  $ and $  r^{M,2}_{ij} \in \mathcal{M}([0,T]\times \T^d)$. 
Moreover, for $ \psi \in C^1([0,T]\times \T^d) $, it holds
\begin{align*}
	\lel E_\infty;\psi  \ril_{\mathcal{M}([0,\tau]\times \T^d), C([0,\tau]\times \T^d)}= \int_0^\tau \lel E_\infty ; \psi \ril_{\mathcal{M}(\T^d), C(\T^d)} \dt .
\end{align*}
The above statement follows from the fact that $ E_\infty \in L^\infty_{\text{weak-(*)}}(0,T;\mathcal{M}(\T^d)) $.
Therefore, we conclude that
\begin{align*}
	\left\vert \lel r^{M,2};\Grad \pmb{\phi}\ril_{\mathcal{M}([0,\tau]\times \T^d), C([0,\tau]\times \T^d)} \right\vert  \leq \frac{1}{4 \epsilon}  \int_{0}^\tau \Vert \pmb{\phi} (t) \Vert_{C^1(\T^d)} E_\infty (\tau ) \dt + \epsilon \Vert \pmb{\phi} \Vert_{C^1([0,\tau]\times \T^d)}  \sigma_\infty[[0,\tau] \times \T^d].
\end{align*}
The total defect is equal to
\[ r^{M}=r^{M,1}+r^{M,2}.\]
We also note that $$ \delta \intTAO{ \varrho_\delta \vw_\delta \cdot \Grad p(\varrho_n) : \Grad \pmb{\phi} }  \rightarrow 0 $$ as $ \delta^\frac{1}{2} \Vert \vw_\delta \Vert_{L^2(0,\tau ;L^2(\T^d))} $ is uniformly bounded. Thus, if we pass to the limit $\delta\to0$ in the momentum equation \eqref{AR2_delta} we get 
\begin{align}\label{mveqn2pf}
	\begin{split}
		&\int_{ \T^d} \sqrt{\vr}(\tau,x) \lel\Nu_{\tau,x}; \rtrw  \ril \cdot \pmb{\phi}(\tau,\cdot) \dx - \int_{ \T^d} \lel \Nu_{0} ;  \sqrt{s} \rtrw  \ril \cdot \pmb{\phi}(0,\cdot) \dx\\
		&= \int_{0}^{\tau} \int_{ \T^d}  \left[ \sqrt{\vr} \lel \Nu_{t,x};  \rtrw   \ril \cdot \partial_{t} \pmb{\phi} +  \lel \Nu_{t,x}; \rtrw \otimes  \rtrw  \ril : \Grad \pmb{\phi}- \lel \Nu_{t,x}; \rtrw \otimes \DQ \ril : \Grad \pmb{\phi} \right] \dx \dt\\
		&  + \lel r^M;\Grad \pmb{\phi}\ril_{\mathcal{M}([0,\tau]\times \T^d), C([0,\tau]\times \T^d)} ,
	\end{split}
\end{align}
for a.e. $\tau \in (0,T)$, for all $ \epsilon>0 $ and $\pmb{\phi}\in C^{1}([0,T]\times {\T^d};\R^d)$, where  $ r^{M}$ satisfies \eqref{mom-def}.

In the energy inequality, there is another nonlinear term $\intTAO{ \varrho_\delta \vw_\delta  \cdot \Grad p(\varrho_\delta) }$. {We define the defect measure related to this term as 
\begin{align*}
	 \mathcal{R}= \overline{\sqrt{\varrho} \vw  \cdot \Grad Q(\varrho)} - \lel \Nu_{t,x}; \rtrw \cdot \DQ \ril,
\end{align*}
where $  \overline{\sqrt{\varrho} \vw  \cdot \Grad Q(\varrho)}  $ stands for the weak-(*) limit of $ \varrho_\delta \vw_\delta  \cdot \Grad p(\varrho_\delta)  $ in $ L^\infty_{\text{weak-(*)}}(0,T; \mathcal{M}(\T^d)) $.} We proceed similarly and for $\epsilon>0$ we have 
\[\left\vert \int_{(0,\tau)\times \T^d}{\text{d} \mathcal{R}}   \right\vert =\vert  \mathcal{R} [[0,\tau] \times \T^d] \vert \leq \frac{C}{4 \epsilon} \int_{0}^\tau E_\infty (t) \dt + \epsilon  \sigma_\infty[[0,\tau] \times \T^d] . \]
This yields
\begin{align*}
			\begin{split}
				&\int_{ \T^d} \lel \Nu_{\tau ,x}; \frac{1}{2}  \levert \rtrw \rivert^2 + E(s) \ril \dx + \int_{0}^{\tau} 	\intO{ \lel \Nu_{t,x}; \vert \DQ \vert^2 \ril} \dt + \mathcal{D}(\tau) \\
				&\leq \int_{ \T^d} \lel \Nu_{0 ,x};  \frac{1}{2} \levert \rtrw \rivert^2 +E(s) \ril \dx + \int_{0}^{\tau} \intO{ \lel \Nu_{t,x}; \rtrw \cdot \DQ \ril }  \dt +\int_{(0,\tau)\times \T^d}{\text{d} \mathcal{R}} .
			\end{split}
		\end{align*}
		The proof of Theorem \ref{Theorem:1} is therefore complete. $\Box$

\section{Relative energy: formal derivation}\label{rel-ent-motivation}
The purpose of this section is to derive the relative energy \eqref{def:rel_ent} and to show that it is sufficient to prove the weak-strong uniqueness property at least at the level of strong solutions.
\subsection{Relative energy inequality}
We assume that $(\vr,\vw)$ is a smooth solution of \eqref{PM}, and we denote by $(\bar\vr,\bar\vw)$ another (smooth) solution to the system \eqref{AR_multi} s.t. $\bar\vr>0$ and $\bar\vu$ is given by \eqref{def:w}. 
Moreover, we introduce the following notation:
\eqh{
\vW=\vw-\bvw,\quad \vU=\vu-\bvu,\quad R=\vr-\bvr.}
First, we verify that the functional given in the introduction by \eqref{def:rel_ent} is a relative energy, i.e., a distance between the solutions $ (\bvr, \bar{\vw}) $ and $ (\vr,\vc{w}) $ controlled by the difference of the corresponding initial data. 
 The first step is to prove the following result:
\begin{lemma}\label{Lemma:3}
Let $(\vr,\vw)$ and $(\bvr,\bvw)$ be smooth solutions to \eqref{AR_multi}, s.t.
\eqh{\bvr>0\quad \mbox{and}\quad \intO{\bvr_0}=\intO{\vr_0}.}
Then for all $t\in [0,T)$ we have
\eqh{
\Dt\intO{\vr\frac{|\vW|^2}{2}}=\intO{\vr\vU\cdot\vW\Grad \bvw}.
}
\end{lemma}
\pf Differentiation with respect to time gives:
\eq{\label{dist_mom}
\Dt\intO{\vr\frac{|\vw-\bar\vw|^2}{2}}=\intO{\pt\vr \frac{|\vw-\bar\vw|^2}{2}}+\intO{\vr(\vw-\bar\vw)\cdot\pt(\vw-\bar\vw)}.
}
The first part can be evaluated using the continuity equation, and the second part from the momentum equations, more precisely
\eqh{
\vr\pt\vw+\vr\vu\Grad\vw=0 \quad \text{ and }\quad \bar\vr\pt\bar\vw+\bar\vr\bar\vu\Grad\bvw=0. }

Therefore, \eqref{dist_mom} is rewritten as
\eq{\label{rel_1}
\Dt\intO{\vr\frac{|\vw-\bvw|^2}{2}}&=-\intO{\Div(\vr\vu)\frac{|\vw-\bvw|^2}{2}}-\intO{\vr(\vw-\bvw)\cdot(\vu\Grad\vw-\bvu\Grad\bvw)}\\
&=\intO{\vr\vu\cdot(\vw-\bvw)\Grad(\vw-\bvw)}-\intO{\vr(\vw-\bvw)\cdot(\vu\Grad\vw-\bvu\Grad\bvw)}\\
&=-\intO{\vr\vu\cdot(\vw-\bvw)\Grad\bvw}+\intO{\vr\bvu\cdot(\vw-\bvw)\Grad\bvw}\\
&=\intO{\vr(\bvu-\vu)\cdot(\vw-\bvw)\Grad \bvw}
}
which concludes the proof. $\Box$

Note that, plugging into this expression formulas $\vu=\vw+\Grad p(\vr)$ and $\bvu=\bvw+\Grad p(\bvr)$, we get that
\eqh{
\Dt\intO{\vr\frac{|\vw-\bvw|^2}{2}}&\leq C_1\intO{\vr|\vw-\bvw|^2|\Grad\bvw|}+C_2\intO{\vr|\Grad p(\vr)-\Grad p(\bvr)|^2|\Grad\bvw|}.
}
The purpose of the next lemma is to control the second term in the above equation. To this purpose we  evaluate the time evolution of the second quantity in \eqref{def:rel_ent}, i.e.
\eqh{
\mathcal{E}(\varrho \mid \bvr)=\intOB{ E(\varrho)-E(\bvr)-E^\prime (\bvr) (\varrho-\bvr) }.
}
\begin{lemma}\label{Lemma:4}
Under the assumptions of Lemma \ref{Lemma:3}, we have
\eq{\label{E:03a}
	\Dt \mathcal{E}(\varrho \mid \bvr) + \intO{ 
		\left\vert \left(\Grad Q(\varrho) -\sqvrf \Grad Q(\bvr) \right) \right \vert^2 } = {\mathfrak{R}}(\vr,\vu,\bvr,\bvu)}
where 		
\eqh{
{\mathfrak{R}}(\vr,\vu,\bvr,\bvu)&=\intO{\sqvr (\vw- \bar{\vw})\cdot\left(\Grad Q(\varrho) -\sqvrf \Grad Q(\bvr) \right)} \\
	&\quad+ \intO{\left(\sqvr- \sqvrb\right)  \bar{\vw}\cdot\left(\Grad Q(\varrho) -\sqvrf \Grad Q(\bvr) \right) }\\
	& \quad+ \intO{ \left(1-\sqvrf \right) \Grad Q(\bar{\varrho})\cdot\left(\Grad Q(\varrho) - \sqvrf \Grad Q(\bvr) \right) }\\
&\quad - \intO{\left( Q(\varrho) - Q(\bvr) -(\vr- \bvr )Q^\prime(\bvr) \right) \Div (\sqvrb \bar{\vu}) }\\
	&\quad- \intO{ \left( \sqvr -\sqvrb - \frac{1}{2 \sqvrb} (\vr - \bvr )\right) \bar{\vu} \cdot \Grad Q(\bar{\varrho})}.
}

\end{lemma}
\pf  We differentiate w.r.t. time to get
\begin{align*}
\Dt \mathcal{E}(\varrho \mid \bvr) =\Dt \intO{E(\varrho)} - \intO{(E^\prime (\bvr) \varrho)_t} + \intO{(\bvr E^\prime(\bvr)-E(\bvr))_t}.
\end{align*}
Using the energy equality \eqref{energy_eq} and mass conservation \eqref{PM1} we obtain 
\begin{align*}
\Dt \mathcal{E}(\varrho \mid \bvr) + \intO{ \varrho 
	\vert \Grad p(\varrho) \vert^2} =&  \intO{\vr\vw\cdot\Grad p(\vr)} - \intO{(\varrho-\bvr) E^{\prime \prime}(\bvr) \bvr_t}\\
&-\intO{\Grad (E^\prime(\bvr) )(\varrho \vc{w} - \varrho \Grad p(\varrho) )}.
\end{align*}
From the fact that $ E^\prime(\vr) = p(\vr) $ the above equality reduces to
\eq{\label{ptoQ}
		&\Dt \mathcal{E}(\varrho \mid \bvr) + \intO{ \varrho 
		 \Grad p(\varrho) (\Grad p(\varrho) -\Grad p(\bvr) )} \\& =  \intO{\vr\vw\cdot(\Grad p(\varrho) -\Grad p(\bvr) )} - \intO{(\varrho-\bvr) E^{\prime \prime}(\bvr) \bvr_t}.
}
We rewrite \eqref{ptoQ} as
\begin{align*}
	&\Dt \mathcal{E}(\varrho \mid \bvr) + \intO{ 
		\Grad Q(\varrho) \left(\Grad Q(\varrho) -\sqvrf \Grad Q(\bvr) \right)} \\& =  \intO{\sqvr \vw\cdot\left(\Grad Q(\varrho) -\sqvrf \Grad Q(\bvr) \right)} - \intO{(\varrho-\bvr) E^{\prime \prime}(\bvr) \bvr_t}.
\end{align*}
Rearranging terms we have 
\eq{\label{E:01}
	&\Dt \mathcal{E}(\varrho \mid \bvr) + \intO{ 
		\left\vert \left(\Grad Q(\varrho) -\sqvrf \Grad Q(\bvr) \right) \right \vert^2 } \\ &=  \intO{\sqvr (\vw- \bar{\vw})\cdot\left(\Grad Q(\varrho) -\sqvrf \Grad Q(\bvr) \right)} \\
	&\quad + \intO{\sqvr \bar{\vw}\cdot\left(\Grad Q(\varrho) -\sqvrf \Grad Q(\bvr) \right)} \\
	&\quad- \intO{ \sqvrf \Grad Q(\bvr) \left(\Grad Q(\varrho) -\sqvrf \Grad Q(\bvr) \right)}\\
	&\quad- \intO{(\varrho-\bvr) E^{\prime \prime}(\bvr) \bvr_t} =\sum_{i=1}^4 I_i.
}
The term $I_1$ already has the form we want. Now, for the term $I_2$ we use \eqref{def:w} to get
\begin{AS}\label{M:ree:1}
	&\sqvr \bar{\vw}\cdot\left(\Grad Q(\varrho) -\sqvrf \Grad Q(\bvr) \right) \\
	&= \left(\sqvr- \sqvrb\right)  \bar{\vw}\cdot\left(\Grad Q(\varrho) -\sqvrf \Grad Q(\bvr) \right) + \sqvrb \bar{\vw}\cdot\left(\Grad Q(\varrho) -\sqvrf \Grad Q(\bvr) \right) \\
	&= \left(\sqvr- \sqvrb\right)  \bar{\vw}\cdot\left(\Grad Q(\varrho) -\sqvrf \Grad Q(\bvr) \right) + {(\sqvrb \bar{\vu}+ \Grad Q(\bar{\varrho}))\cdot\left(\Grad Q(\varrho) - \sqvrf \Grad Q(\bvr) \right) }.
\end{AS}
We rewrite the last term in the above equation in the following way
\begin{align}\label{M:ree:2}
\begin{split}
	& \Grad Q(\bar{\varrho})\cdot\left(\Grad Q(\varrho) - \sqvrf \Grad Q(\bvr) \right) \\
	 &= \sqvrf  \Grad Q(\bar{\varrho})\cdot\left(\Grad Q(\varrho) - \sqvrf \Grad Q(\bvr) \right) + \left(1-\sqvrf \right) \Grad Q(\bar{\varrho})\cdot\left(\Grad Q(\varrho) - \sqvrf \Grad Q(\bvr) \right), 
\end{split}
\end{align}
so that the first part cancels out $I_3$, and 
\begin{align*}
	& \sqvrb \bar{\vu}  \left(\Grad Q(\varrho) - \sqvrf \Grad Q(\bvr) \right)=  \sqvrb \bar{\vu}  \left(\Grad Q(\varrho) - \Grad Q(\bvr) \right) - \left(\sqvrf-1\right) \sqvrb \bar{\vu}   \Grad Q(\bvr).
\end{align*}
The first term can be integrated by parts, thus
\begin{align}\label{M:ree:3}
\begin{split}
	\intO{\sqvrb \bar{\vu} \left(\Grad Q(\varrho) - \Grad Q(\bvr) \right) }= - \intO{\Div(\sqvrb \bar{\vu}) \left( Q(\varrho) -  Q(\bvr) \right)}.
\end{split}
\end{align}
Plugging these calculations into \eqref{E:01} we obtain
\eq{\label{E:02}
	&\Dt \mathcal{E}(\varrho \mid \bvr) + \intO{ 
		\left\vert \left(\Grad Q(\varrho) -\sqvrf \Grad Q(\bvr) \right) \right \vert^2 } \\
	&=  \intO{\sqvr (\vw- \bar{\vw})\cdot\left(\Grad Q(\varrho) -\sqvrf \Grad Q(\bvr) \right)} \\
	&\quad+ \intO{\left(\sqvr- \sqvrb\right)  \bar{\vw}\cdot\left(\Grad Q(\varrho) -\sqvrf \Grad Q(\bvr) \right) }\\
	& \quad+ \intO{ \left(1-\sqvrf \right) \Grad Q(\bar{\varrho})\cdot\left(\Grad Q(\varrho) - \sqvrf \Grad Q(\bvr) \right) }\\
	& \quad- \intO{\Div(\sqvrb \bar{\vu}) \left( Q(\varrho) -  Q(\bvr) \right) } - \intO{\left(\sqvrf-1\right) \sqvrb \bar{\vu}   \Grad Q(\bvr)}\\
	&\quad-\intO{(\varrho-\bvr) E^{\prime \prime}(\bvr) \bvr_t}.
}
For the last term in the above equation, we use the continuity equation for the strong solutions to get 
\begin{AS}\label{M:ree:4}
      -\intO{(\varrho-\bvr) E^{\prime \prime}(\bvr) \bvr_t}= \intO{(\varrho-\bvr) E^{\prime \prime}(\bvr) \Div(\bvr \bar{\vu})}. 
\end{AS}
Using the relation 
\[ \frac{1}{\sqvrb} \Div (\bar{\varrho}\bar{\vu}) =  \Div (\sqvrb \bar{\vu}) + \bar{\vu} \cdot \Grad \sqvrb,\] 
along with relations $ E^\prime (\vr)=p(\vr) $ and $Q^\prime(\varrho) =\sqrt{\varrho} p^\prime(\varrho),
$
we obtain
\begin{AS}\label{M:ree:5}
   	\intO{(\varrho-\bvr) E^{\prime \prime}(\bvr) \Div(\bvr \bar{\vu})}&=	\intO{(\varrho-\bvr) Q^{\prime }(\bvr)  \frac{1}{\sqvrb} \Div(\bvr \bar{\vu}) }\\
		&	= \intO{(\varrho-\bvr) Q^{\prime }(\bvr) (\Div (\sqvrb \bar{\vu}) + \bar{\vu} \cdot \Grad \sqvrb) }. 
\end{AS}
At last, combining the last three terms, we have
\eq{\label{E:03}
	&\Dt \mathcal{E}(\varrho \mid \bvr) + \intO{ 
		\left\vert \left(\Grad Q(\varrho) -\sqvrf \Grad Q(\bvr) \right) \right \vert^2 } \\
	&\leq  \intO{\sqvr (\vw- \bar{\vw})\cdot\left(\Grad Q(\varrho) -\sqvrf \Grad Q(\bvr) \right)} \\
	&\quad+ \intO{\left(\sqvr- \sqvrb\right)  \bar{\vw}\cdot\left(\Grad Q(\varrho) -\sqvrf \Grad Q(\bvr) \right) }\\
	& \quad+ \intO{ \left(1-\sqvrf \right) \Grad Q(\bar{\varrho})\cdot\left(\Grad Q(\varrho) - \sqvrf \Grad Q(\bvr) \right) }\\
&\quad - \intO{\left( Q(\varrho) - Q(\bvr) -(\vr- \bvr )Q^\prime(\bvr) \right) \Div (\sqvrb \bar{\vu}) }\\
	&\quad- \intO{ \left( \sqvr -\sqvrb - \frac{1}{2 \sqvrb} (\vr - \bvr )\right) \bar{\vu} \cdot \Grad Q(\bar{\varrho})}
}
which concludes the proof. $\Box$

Finally, plugging  the outcomes of Lemmas \ref{Lemma:3} and \ref{Lemma:4} into the definition of relative energy \eqref{def:rel_ent}, we conclude the subsection with the following lemma.
\begin{lemma}\label{Lemma:13}
	Suppose $ (\vr,\vw) $ and $ (\bvr,\bvw) $ are two sets of strong solutions for the system \eqref{PM}. Then we have the following inequality:
	\begin{AS}\label{E:04}
	        	&\Dt \mathcal{E}(\varrho, \vc{w} \mid \bvr , \bar{\vw})  + \intO{ 
		\left\vert \left(\Grad Q(\varrho) -\sqvrf \Grad Q(\bvr) \right) \right \vert^2 } \\
	& \leq-\intO{\vr(\vw -\bvw)\cdot(\vw-\bvw)\Grad \bvw}\\
	& + \intO{ \sqvr (\vw- \bvw) \cdot \left(\Grad Q(\varrho) -\sqvrf \Grad Q(\bvr) \right) \Grad \bvw }\\
	& +  \intO{\sqvr (\vw- \bar{\vw})\cdot\left(\Grad Q(\varrho) -\sqvrf \Grad Q(\bvr) \right)} \\
	&+ \intO{\left(\sqvr- \sqvrb\right)  \bar{\vw}\cdot\left(\Grad Q(\varrho) -\sqvrf \Grad Q(\bvr) \right) }\\
	& + \intO{ \left(1-\sqvrf \right) \Grad Q(\bar{\varrho})\cdot\left(\Grad Q(\varrho) - \sqvrf \Grad Q(\bvr) \right) }\\
	& - \intO{\left( Q(\varrho) - Q(\bvr) -(\vr- \bvr )Q^\prime(\bvr) \right) \Div (\sqvrb \bar{\vu}) }\\
	&- \intO{ \left( \sqvr -\sqvrb - \frac{1}{2 \sqvrb} (\vr - \bvr )\right) \bar{\vu} \cdot \Grad Q(\bar{\varrho})}= \sum_{i=1}^{7}\mathcal{T}_i.
	\end{AS}
\end{lemma}
\subsection{Relative energy estimate}
Our goal is to show that the terms in \eqref{E:04} are either controlled by the relative energy or can be absorbed in the left hand side.
To do this, we first give a lemma which is a direct consequence of the coercivity of the relative energy functional.
\begin{lemma}
	Suppose $E^\prime (\vr) = p(\vr)$ and $\bvr$ lies in a compact subset of $(0,\infty)$, then for some $ 0<r_1<r_2$ we have
	\begin{align*}
		E(\varrho)-E(\bvr)-E^{\prime}(\bvr)(\varrho -\bvr) \geq  
		\begin{cases}
			&c(\bvr) (\varrho -r)^2 \text{ for } r_1 \leq \varrho \leq r_2,\\
			&c(\bvr) (1+ \varrho^{\gamma+1} ) \text{ otherwise }\\
		\end{cases},
	\end{align*}
	where, $c(\bvr)$ is uniformly bounded in $ [0,\infty)$ for $\bvr$ belonging to compact subsets of $(0,\infty)$.
\end{lemma}
In our case, we consider $r_1,r_2>0$ such that they satisfy 
$$r_1 < \frac{\inf\limits_{(x,t)\in (0,T)\times \T^d} \bvr(x,t)}{2},\qquad r_2 > 2 \times {\sup\limits_{(x,t)\in (0,T)\times \T^d} \bvr(x,t)} $$
and
$$1 + \varrho^\gamma \geq \max\{ \varrho^{\gamma + \frac{1}{2}}, \varrho^2 \},\text{ for all } \varrho \geq r_2 .$$
Now, we estimate each term in the  r.h.s of \eqref{E:04}. First we notice that
\begin{align*}
	\vert \mathcal{T}_1 \vert \leq \Vert \Grad \bvw \Vert_{L^\infty(\T^d)} \intO{\varrho \vert \vw -\bvw \vert^2 }.
\end{align*}
Then for $ \mathcal{T}_2 $ and $ \mathcal{T}_3 $ we obtain 
\begin{align*}
		\vert \mathcal{T}_2 + \mathcal{T}_3 \vert  \leq C_\eta \left(1+\Vert \Grad \bvw \Vert_{L^\infty(\T^d)}\right)^2 \intO{\varrho \vert \vw -\bvw \vert^2 }+ \eta \intO{ 
			\left\vert \left(\Grad Q(\varrho) -\sqvrf \Grad Q(\bvr) \right) \right \vert^2 } ,
\end{align*}
where $ \eta >0 $. Since $ \bvr $ is a strong solution with $ 0< \inf_{(0,T)\times \T^d} \bvr  = \underline{\vr}$, we notice that 
\begin{align}\label{M:ws:1}
	\vert \sqvr- \sqvrb \vert \leq \frac{\vert \vr - \bvr \vert }{ \sqvr+ \sqvrb }  \leq  \frac{\vert \vr - \bvr \vert }{  \sqrt{\underline{\vr}}  } .
\end{align}
Thus, for $ \mathcal{T}_4 $ and $ \mathcal{T}_5 $ we have 
\begin{align*}
		\vert \mathcal{T}_4 + \mathcal{T}_5 \vert  \leq& \frac{C_\eta }{\underline{\vr}}\left(\Vert Q(\bvr) \Vert_{L^\infty(\T^d)}+\Vert \bvw \Vert_{L^\infty(\T^d)}\right)^2 \intO{ \vert \vr -\bvr \vert^2 }\\
		&+ \eta \intO{ 
		\left\vert \left(\Grad Q(\varrho) -\sqvrf \Grad Q(\bvr) \right) \right \vert^2 }. 
\end{align*}
Since $ p(\vr) = \varrho^\gamma  $ with $ \gamma\geq 1 $ we have $ E \approx \varrho^{\gamma+1} $ thus $(\vr -\bvr)^2$ is controlled by $ E(\vr\vert \bvr) $. Furthermore, using $ Q \approx \varrho^{\gamma+ \frac{1}{2}} $, we have 
\begin{align*}
	\vert \mathcal{T}_6  \vert  \leq \Vert \sqvrb \bar{\vu} \Vert_{L^\infty(\T^d)} E(\vr\vert \bvr).
	\end{align*}
For the term $ \mathcal{T}_7 $, first we note that 
\begin{align*}
	 \left( \sqvr -\sqvrb - \frac{1}{2 \sqvrb} (\vr - \bvr )\right) &= (\varrho - \bvr) \left( \frac{1}{\sqvr + \sqvrb} - \frac{1}{2 \sqvrb} \right)=(\varrho - \bvr) \frac{\sqvrb-\sqvr}{2(\sqvr + \sqvrb) \sqvrb}\\
	 &=- \frac{(\varrho - \bvr)^2}{2(\sqvr + \sqvrb)^2 \sqvrb}.
\end{align*}
Hence we have the estimate
\begin{align}\label{M:ws:2}
	\left\vert  \left( \sqvr -\sqvrb - \frac{1}{2 \sqvrb} (\vr - \bvr )\right)  \right\vert \leq \frac{\vert \vr -\bvr \vert^2}{2 \underline{\vr}^{\frac{3}{2}}},
\end{align}
and so, as a consequence
\begin{align*}
	\vert \mathcal{T}_7  \vert  \leq \frac{1}{\underline{\vr}^{\frac{3}{2}}} E(\vr\vert \bvr).
\end{align*}
Ultimately,  choosing $ \delta $ suitably small we deduce 
\begin{align*}
	&\Dt  \mathcal{E}(\varrho, \vc{w} \mid \bvr , \bar{\vw}) + \frac{1}{2}\intO{ 
		\left\vert \left(\Grad Q(\varrho) -\sqvrf \Grad Q(\bvr) \right) \right \vert^2 } \leq C(\bvr,\bvw) \mathcal{E}(\varrho, \vc{w} \mid \bvr , \bar{\vw})  .
\end{align*}
We use Gr\" onwall's lemma to conclude $ \vr = \bvr $ and $ \vw= \bvw $ if they share same initial data. 


\section{Weak(Measure-valued)-Strong uniqueness}\label{Sec:6}
In Section \ref{rel-ent-motivation}, we gave a formal derivation of the relative energy inequality for the system \eqref{PM}. It is clear that for the derivation we need certain regularity for  $ \vw $, which is missing in Definition \ref{MVdef}.\par A suitable relative energy in  terms of the Young measure is given by
\eq{\label{rel-ent-mv}
	\mathcal{E}_{mv}(\tau)&=	\mathcal{E}_{mv}(\Nu, \mathcal{D} \vert \bvr,\bvw)(\tau) \\
	&:= \intO{\lel \Nu_{\tau ,x};\frac{1}{2} \levert\rtrw - \sqrt{s} \bvw \rivert^2 +\left( E(s)- E(\bvr)-E^\prime(\bvr)(s-\bvr)\right) \ril (\tau )},
}
where $(\Nu,\mathcal{D}) $ is the measure-valued solution and $ (\bvr,\bvw) \in C^2([0,T]\times\T^d) \times C^2([0,T]\times\T^d; \R^d) $ with $\bvr >0 $.
\subsection{Relative energy inequality for measure-valued solutions}
For the measure-valued solution we have the following analogue of Lemma \ref{Lemma:13}.
\begin{lemma}
	Let $ \{\Nu, \mathcal{D}\} $ be a measure-valued solution and $ \{\bvr,\bvw\} $ be a strong solution of the system \eqref{PM} with $ \bvr >0 $ in $ (0,T)\times \T^d $. Then we have 
	\begin{AS}\label{MS:re0}
			\mathcal{E}_{mv}(\tau) + & 	\intTAO{ \lel \Nu_{t,x};  \levert\DQ-\sqvrft\Grad Q(\bvr) \rivert^2 \ril}  + \mathcal{D}(\tau) \\
			&\leq \mathcal{E}_{mv}(0) + \mathfrak{R}_{mv}(\tau) ,  
		\end{AS} 
		where the reminder is of the form
		\begin{AS}\label{MS:re1}
		   \mathfrak{R}_{mv} =&     \intTAO{ \lel \Nu_{t,x}; \left(\rtrw- \sqrt{s}\bvw\right)  \cdot \left(\DQ-\sqvrft\Grad Q(\bvr) \right)  \ril }   \\
			&- \intTAO{\lel \Nu_{t,x}; \left(\rtrw- \sqrt{s}\bvw\right)  \otimes \left(\rtrw- \sqrt{s}\bvw\right)   \ril : \Grad \bvw}\\
			&+ \intTAO{\lel \Nu_{t,x}; \left(\rtrw- \sqrt{s}\bvw\right) \otimes \left(\DQ-\sqvrft\Grad Q(\bvr) \right)   \ril :  \Grad \bvw} \\
			&+ \intTAO{\lel \Nu_{t,x}; \left(\sqrt{s}- \sqvrb \right) \left(\DQ -\sqvrft \Grad Q(\bvr) \right) \ril  \cdot \bar{\vw} }\\
			& + \intTAO{ \lel \Nu_{t,x}; \left(1-\sqvrft \right) \Grad Q(\bar{\varrho})\cdot\left(\DQ - \sqvrft \Grad Q(\bvr) \right) \ril }\\
			& - \intTAO{\lel \Nu_{t,x};\ \left( Q( {s}) - Q(\bvr) -( {s}- \bvr )Q^\prime(\bvr) \right)\ril  \Div (\sqvrb \bar{\vu}) }\\
			&- \intTAO{\lel \Nu_{t,x};\ \left( \sqrt{s}-\sqvrb - \frac{1}{2 \sqvrb} ( {s}- \bvr )\right) \ril \bar{\vu} \cdot \Grad Q(\bar{\varrho})}\\
			&   +   \lel r^M;\Grad \pmb{\phi}\ril_{\mathcal{M}([0,\tau]\times \T^d), C([0,\tau]\times \T^d)}   +\int_{(0,\tau)\times \T^d}{\text{d} \mathcal{R}}
		\end{AS}
  with $ \bar{\vu}= \bar{\vw}- \Grad p(\bvr)$.
\end{lemma}
\pf
 At first, a simple calculation leads us to the following identity:
\begin{align}\label{relm0}
	\begin{split}
	\mathcal{E}_{mv}(\tau) = &\intO{\lel \Nu_{\tau ,x};  \frac{1}{2} \levert \rtrw \rivert^2 + E(s) \ril } - \intO{\sqrt{\vr}(\tau, x)\lel \Nu_{\tau ,x};  \rtrw \ril \cdot \bvw(\tau ,x)} \\
	& +\intO{\vr(\tau, x) \left( \frac{\vert \bvw \vert^2}{2}- E^\prime(\bvr) \right)(\tau,x)}+\intO{(\bvr E^\prime(\bvr)-E(\bvr))(\tau ,x)}= \sum_{i=1}^{4} \mathcal{J}_i^{(1)}.
\end{split}
\end{align}
In the above identity, the first term $\mathcal{J}_1^{(1)} $ will be replaced by the energy inequality from the Definition \ref{MVdef}. For $ \mathcal{J}_2^{(1)} $ and $ \mathcal{J}_3^{(1)}$, we use \eqref{mv eqn 1} and \eqref{mv eqn 2}, respectively. Finally, we rewrite the term $ \mathcal{J}_4^{(1)}$ as
\begin{align*}
	\mathcal{J}_4^{(1)} & =\int_0^\tau \intO{\partial_{t} (\bvr E^\prime(\bvr)-E(\bvr)) }\dt + \intO{(\bvr E^\prime(\bvr)-E(\bvr))(0)}.
\end{align*}
After substituting for $ \mathcal{J}_i^{(1)} $'s into \eqref{relm0}, we obtain the inequality
\begin{align*}
	&\mathcal{E}_{mv}(\tau) + \intTAO{ \lel \Nu_{t,x};  \DQ \left(\DQ-\sqvrft\Grad Q(\bvr) \right) \ril} + \mathcal{D}(\tau) \\
	&\leq \mathcal{E}_{mv}(0) +  \intTAO{ \lel \Nu_{t,x}; \rtrw \cdot \left(\DQ-\sqvrft\Grad Q(\bvr) \right)  \ril }   \\
	&- \intTAO{\lel \Nu_{t,x}; \sqrt{s} \left(\rtrw- \sqrt{s}\bvw\right)\ril \bvw_t}\\
	&-\intTAO{\lel \Nu_{t,x}; \left(\rtrw- \sqrt{s}\bvw\right) \otimes \rtrw \ril :  \Grad \bvw} \\
	&+ \intTAO{\lel \Nu_{t,x}; \sqrt{s} \left(\rtrw- \sqrt{s}\bvw\right) \otimes \DQ \ril :  \Grad \bvw} \\
	&- \intTAO{\lel\Nu_{t,x}; s- \bvr \ril \left( E^\prime (\bvr) \right)_t}\\
	& +   \lel r^M;\Grad \pmb{\phi}\ril_{\mathcal{M}([0,\tau]\times \T^d), C([0,\tau]\times \T^d)}  +\int_{(0,\tau)\times \T^d}{\text{d} \mathcal{R}} .
\end{align*}
Next, since $ (\bvr,\bvw) $ is a strong solution of the system \eqref{PM}, the above inequality reduces to
\begin{AS}\label{re2}
    	\mathcal{E}_{mv}(\tau) + & 	\intTAO{ \lel \Nu_{t,x};  \DQ \left(\DQ-\sqvrft\Grad Q(\bvr) \right) \ril}  + \mathcal{D}(\tau) \\
	&\leq \mathcal{E}_{mv}(0) +   \intTAO{ \lel \Nu_{t,x}; \left(\rtrw- \sqrt{s}\bvw\right)  \cdot \left(\DQ-\sqvrft\Grad Q(\bvr) \right)  \ril }   \\
	&- \intTAO{\lel \Nu_{t,x}; \left(\rtrw- \sqrt{s}\bvw\right)  \otimes \left(\rtrw- \sqrt{s}\bvw\right)   \ril : \Grad \bvw}\\
	&+ \intTAO{\lel \Nu_{t,x}; \left(\rtrw- \sqrt{s}\bvw\right) \otimes \left(\DQ-\sqvrft\Grad Q(\bvr) \right)   \ril :  \Grad \bvw} \\
	&+ \intTAO{ \lel \Nu_{t,x}; \sqrt{s} \bvw \cdot \left(\DQ-\sqvrft\Grad Q(\bvr) \right)  \ril }  \\
	&- \intTAO{\lel\Nu_{t,x}; s- \bvr \ril \left( E^\prime (\bvr) \right)_t}\\
	&   +   \lel r^M;\Grad \pmb{\phi}\ril_{\mathcal{M}([0,\tau]\times \T^d), C([0,\tau]\times \T^d)}   +\int_{(0,\tau)\times \T^d}{\text{d} \mathcal{R}}= \sum_{i=0}^{7}\mathcal{J}^{(2)}_i.   
\end{AS}
For $ \mathcal{J}^{(2)}_4 $, we can perform rearrangement of the terms as we did in \eqref{M:ree:1} and \eqref{M:ree:2}.
Next, we use the $L^2$ integrability of the term $ \Grad Q(\varrho) $ to obtain
\begin{align}\label{MS:1}
	\intO{\sqvrb \bar{\vu} \left(\lel \Nu_{t,x}; \DQ \ril- \Grad Q(\bvr) \right) }= - \intO{\Div(\sqvrb \bar{\vu}) \left( \lel \Nu_{t,x}; Q(s)\ril  -  Q(\bvr) \right) }.
\end{align}
Also, for the term $ \mathcal{J}^{(2)}_5 $, we follow similar calculation as in \eqref{M:ree:4} and \eqref{M:ree:5}. 
Finally, we derive the relative energy inequality as stated in \eqref{MS:re0} with remainder term $ \mathfrak{R}_{mv}$ as in \eqref{MS:re1}. 

We split 
\begin{align*}
\mathfrak{R}_{mv}= \sum_{i=0}^{9} \mathcal{J}^{(3)}_i. 
\end{align*}
The expressions for the consecutive parts of the reminder can be pretty much copied from the formal level. Indeed, most of the identities in \eqref{M:ree:1}, \eqref{M:ree:2}, \eqref{M:ree:4} and \eqref{M:ree:5} are either rearrangements or algebraic calculations easily justified on the measure-valued level. The only difficulty is related to the analogue of the identity \eqref{M:ree:3}, but it is justified in \eqref{MS:1}. $\Box$

\medskip
\subsection{Proof of Theorem \ref{Theorem:2}}
To conclude the proof of Theorem \ref{Theorem:2}, we need to establish that the terms in $\mathfrak{R}_{mv}$ are either controlled by the relative energy or can be absorbed in the left hand side of \eqref{MS:re0}.
Since measure-valued and strong solution share same initial data thus we have 
\begin{align}\nonumber
	\mathcal{J}^{(3)}_0  =0.
\end{align}
We also get
\begin{align*}
	\vert \mathcal{J}^{(3)}_2   \vert \leq \Vert \Grad \bvw \Vert_{L^\infty((0,T)\times \T^d)} \intTAO{ \lel \Nu_{t,x}; \levert\rtrw - \sqrt{s} \bvw \rivert^2 \ril } \leq \Vert \Grad \bvw \Vert_{L^\infty((0,T)\times \T^d)} \int_0^\tau \mathcal{E}_{mv}(t) \dt .
\end{align*}
Then, for $ \mathcal{J}^{(3)}_1 $ and $ \mathcal{J}^{(3)}_3 $, we obtain 
\begin{align*}
	\vert \mathcal{J}^{(3)}_1   + \mathcal{J}^{(3)}_3\vert  &\leq C_\eta \left(1+\Vert \Grad \bvw \Vert_{L^\infty(\T^d)}\right)^2 \intTAO{\lel \Nu_{t,x}; \levert\rtrw - \sqrt{s} \bvw \rivert^2  \ril}\\
	&+ \eta \intTAO{ \lel \Nu_{t,x};
		\left\vert \left(\DQ-\sqvrft \Grad Q(\bvr) \right) \right \vert^2 \ril} ,
\end{align*}
for some $ \eta >0$. Now, for the term $\mathcal{J}^{(3)}_4   $ and $ \mathcal{J}^{(3)}_5 $ we follow the computation from \eqref{M:ws:1} and obtain 
\begin{align*}
	\vert \mathcal{J}^{(3)}_4 + \mathcal{J}^{(3)}_5  \vert  \leq& \frac{C_\eta }{\underline{\vr}}\left(\Vert Q(\bvr) \Vert_{L^\infty(\T^d)}+\Vert \bvw \Vert_{L^\infty(\T^d)}\right)^2 \intTAO{ \lel \Nu_{t,x}; \vert \vr -\bvr \vert^2 \ril }\\
	&+ \eta \intTAO{ \lel \Nu_{t,x};
		\left\vert \left(\DQ-\sqvrft \Grad Q(\bvr) \right) \right \vert^2 \ril}, 
\end{align*}
where $ \eta >0 $ and $ 0< \inf_{(0,T)\times \T^d} \bvr  = \underline{\vr}$.
Since $ p \approx \varrho^\gamma  $ with $ \gamma\geq 1 $, we have $ E \approx \varrho^{\gamma+1} $ thus $(\vr -\bvr)^2$ is controlled by $ \mathcal{E}_{mv} $.
Moreover, $ Q \approx \varrho^{\gamma+ \frac{1}{2}} $ leads to 
\begin{align*}
\vert \mathcal{J}^{(3)}_6	\vert  \leq \Vert \sqvrb \bar{\vu} \Vert_{L^\infty((0,T)\times \T^d)} \int_0^\tau \mathcal{E}_{mv}(t) \dt .
\end{align*}
Following the computation from \eqref{M:ws:2}, we estimate
\begin{align*}
	\vert \mathcal{J}^{(3)}_7 \vert  \leq C(\bvr,\bar{\vu})  \int_0^\tau \mathcal{E}_{mv}(t) \dt .
\end{align*}
From the properties of the defect measures \eqref{mom-def} and \eqref{En:def}, for any $ \epsilon >0$ we get
\begin{align*}
	\vert \mathcal{J}^{(3)}_8 +  \mathcal{J}^{(3)}_9 \vert \leq C(\bvw)\frac{1}{4\epsilon} \int_{0}^{\tau} \mathcal{D}(t) \dt  + \epsilon  \left(\levertl \Grad  \bvw \rivertl_{C([0,T]\times \T^d}) +1\right) \mathcal{D}(\tau) .
\end{align*}
Choosing $ \eta $ and $ \epsilon $ sufficiently small, we ultimately deduce
\begin{align*}
	\mathcal{E}_{mv}(\tau) + &\frac{1}{2} 	\intTAO{ \lel \Nu_{t,x};  \left\vert \left(\DQ-\sqvrft\Grad Q(\bvr) \right)\right\vert^2 \ril}  + \frac{1}{2} \mathcal{D}(\tau) \\
	&\leq \mathcal{E}_{mv}(0) +  C(\bvr,\bvw)  \int_0^\tau \mathcal{E}_{mv}(t) \dt + C(\bvw) \int_{0}^{\tau} \mathcal{D}(t) \dt .
\end{align*}
Thus, Gr\" onwall's argument gives the desired result and completes the proof of Theorem \ref{Theorem:2}. $\Box$

\bigskip

{\bf Acknowledgements.} The research of N.C. and of E.Z. leading to these results has received funding from the EPSRC Early Career Fellowship no. EP/V000586/1.
The work of P. G. was supported by National Science Centre (Poland), agreement no 2017/27/B/ST1/01569.
This work was also partially supported by the Simons Foundation Award No 663281 granted to the Institute of Mathematics of the Polish Academy of Sciences for the years 2021-2023.

%
%
%
%
%


\begin{thebibliography}{99}
\bibitem{AFN21}
A.~Abbatiello, E.~Feireisl, and A.~Novotn\'y.
\newblock {Generalized solutions to models of compressible viscous fluids.}
\newblock {\em Discrete \& Continuous Dynamical Systems},  41(1): 1--28, 2021.

\bibitem{AwKlar}
A.~Aw, A.~Klar, M.~Rascle, and T.~Materne.
\newblock Derivation of continuum traffic flow models from microscopic follow-the-leader models.
\newblock {\em SIAM J. Math. Anal.}, 63(1):259--278, 2002.

	
	\bibitem{AR}
	A. Aw and M. Rascle. 
	\newblock Resurrection of second order models of traffic flow.
	\newblock {\em SIAM J. Appl. Math.}, 60:916--938, 2000.
\bibitem{BM}
\newblock{J. M. Ball and F. Murat}
\newblock{Remarks on Chacon's biting lemma.}
\newblock{\em Proc. Amer. Math. Soc. 107 (1989), no. 3, 655–663.}

\bibitem{Basaric}
 \newblock{D. Basari\' c}
 \newblock{ Vanishing viscosity limit for the compressible Navier–Stokes system via measure-valued solutions.}
 \newblock{ \em NoDEA Nonlinear Differential Equations Appl. 27 (2020)}
\bibitem{BrFeHo}	
D. Breit, E. Feireisl, and M. Hofmanov\'a.
\newblock Generalized solutions to models of inviscid fluids.
\newblock {\em Disc. Cont. Dynam. Sys. - B}, 25(10):3831--3841, 2020.
	
	
	\bibitem{Brenier}
	Y. Brenier, C. De Lellis, and L. Sz\'ekelyhidi, Jr. 
	\newblock Weak-strong uniqueness for measure-valued solutions. \newblock {\em Comm. Math. Phys.}, 305(2):351--361, 2011.
	
	\bibitem{B1}
	H. Brenner. 
	\newblock Navier-Stokes revisited. 
	\newblock {\emph Phys. A}, 349(1-2):60–132, 2005.
	
	\bibitem{B2}
	H. Brenner. 
	\newblock Fluid mechanics revisited. 
	\newblock{\emph Phys. A}, 349:190–224, 2006.
	
	
	\bibitem{BD}
D.~Bresch and B.~Desjardins.
\newblock Existence of global weak solutions for a 2d viscous shallow water equations and convergence to the quasi-geostrophic model.
\newblock {\em Commun. Math. Phys.}, 238(1):211--223, 2003.

\bibitem{BrDe}
D.~Bresch and B.~Desjardins.
\newblock On the existence of global weak solutions to the {N}avier-{S}tokes
  equations for viscous compressible and heat conducting fluids.
\newblock {\em J. Math. Pures Appl. }, 87(1):57--90, 2007.

\bibitem{BrDeZa}
D.~Bresch, B.~Desjardins, and E.~Zatorska.
\newblock Two-velocity hydrodynamics in fluid mechanics, part {II}: {E}xistence of global $\kappa$-entropy solutions to the compressible {N}avier--{S}tokes systems with degenerate viscosities.
\newblock {\em J. Math. Pures Appl.}, 104(4):801 -- 836, 2015.

\bibitem{BrGiZa}
D.~Bresch, V.~Giovangigli, and E.~Zatorska.
\newblock Two-velocity hydrodynamics in fluid mechanics, part {I}: {W}ell posedness for zero {M}ach number systems.
\newblock {\em J. Math. Pures Appl.}, 104(4):762 -- 800, 2015.


\bibitem{BrNoVi2017}
D.~Bresch, P.~Noble, and J.-P. Vila.
\newblock Relative entropy for compressible {N}avier-{S}tokes equations with density dependent viscosities and various applications.
\newblock {\em ESAIM: ProcS}, 58:40--57, 2017.
	
		\bibitem{BFN}
	J. B{\v r}ezina, E. Feireisl, and A. Novotn\'y. 
	\newblock Stability of strong solutions to the Navier-Stokes-Fourier system. 
	\newblock {\em SIAM Journal on Mathematical Analysis}, 52(2), 1761--1785, 2020. 
	
	\bibitem{CDGS}
	J. A. Carrillo, T. Debiec, P. Gwiazda and A. \'Swierczewska-Gwiazda.
	\newblock Dissipative measure-valued solutions to the Euler-Poisson equation.
	\newblock{\em arXiv:2109.07536v1}.
	
	\bibitem{CaWrZa2}
J.~A. Carrillo, A.~Wr{\'o}blewska-Kami{\'n}ska, and E.~Zatorska.
\newblock Pressureless {E}uler with nonlocal interactions as a singular limit of degenerate {N}avier-{S}tokes system.
\newblock {\em J. Math. Anal. Appl.}, 492(1), 2020.
	
	
	
	\bibitem{Ch20}
	N. Chaudhuri.
	\newblock On weak (measure-valued)-strong uniqueness for compressible Navier-Stokes system with non-monotone pressure law.  
	\newblock {\em J. Math. Fluid Mech.} 22, 17 (2020).
	
	\bibitem{ChTz}
		C. Christoforou and A. E. Tzavaras. 
		\newblock Relative entropy for hyperbolic-parabolic systems and application to the constitutive theory of thermoviscoelasticity. 
		{\em Arch. Ration. Mech. Anal.}, 229(1):1–52, 2018.
		
			\bibitem{D1979}
			C. Dafermos. 
	\newblock {The second law of thermodynamics and stability}. 
	\newblock {\em Arch. Rational Mech. Anal.}, 70:167–-179, 1979.
	
	\bibitem{Daganzo}
	C. Daganzo. 
	\newblock Requiem for second order fluid approximations of traffic flow. 
	\newblock{\emph Transp. Res. B} 29B: 277-–286, 1995.
	
\bibitem{kis}
T.~Do, A.~Kiselev, L.~Ryzhik, and C.~Tan.
\newblock Global regularity for the fractional {E}uler alignment system.
\newblock {\em Arch. Ration. Mech. Anal.}, 228(1):1--37, 2018.
	
\bibitem{FGSW}
E.~Feireisl, P.~Gwiazda, A.~{\'S}wierczewska--Gwiazda and E.~{Wiedemann}.
\newblock { Dissipative measure-valued solutions to the compressible Navier-Stokes system},
\newblock { \em Calc. Var. Partial Differential Equations}, 55(6):Art. 141, 20, 2016.
	
	

\bibitem{FJN2012}
E. Feireisl, B. J. Jin, and A. Novotn\'y. 
\newblock   {Relative entropies, suitable weak solutions and weak-strong uniqueness for the compressible Navier-Stokes system. }
\newblock {\em J. Math. Fluid Mech.}, 14(4):717--730, 2012.


\bibitem{FL}
E. Feireisl and M. Luk\'a{\r c}ov\'a-Medvid’ov\'a. 
\newblock { Convergence of a mixed finite element–finite volume scheme for the isentropic Navier–Stokes system via dissipative measure-valued solutions.} 
\newblock{\em Foundations of Computational Mathematics}, 18(3):703--730, 2018.

\bibitem{FLM}
E. Feireisl, M. Luk\'a{\r c}ov\'a-Medvid’ov\'a, and H. Mizerov\'a. \newblock { Convergence of finite volume schemes for the Euler equations via dissipative measure-valued solutions.} 
\newblock{\em  Foundations of Computational Mathematics}, 20(4):923-- 966, 2020.
	
\bibitem{FLMS}
E. Feireisl, M. Luk\'acov\'a-Medvid’ov\'a, H. Mizerov\'a, and B. She.
\newblock{\em Numerical Analysis of Compressible Fluid Flows.}
\newblock Springer, 2021.
	

	
\bibitem{FN}
E.~Feireisl and A.~Novotn{{\'y}}.
\newblock {\em Singular limits in thermodynamics of viscous fluids}.
\newblock Advances in Mathematical Fluid Mechanics. Birkh{\"a}user Verlag, Basel, 2009.

\bibitem{FN2021}
E. Feireisl and A. Novotn\'y
\newblock{Weak–strong uniqueness property for models of compressible viscous fluids near vacuum,}
\newblock{\em Nonlinearity, 2021}
	
\bibitem{FNS2011} E. Feireisl, A. Novotn\'y and Y. Sun. 
\newblock  {Suitable weak solutions to the Navier- Stokes equations of compressible viscous fluids}. 
\newblock{\em Indiana Univ. Math. J.}, 60(2):611–631, 2011.

\bibitem{FV}
E.~Feireisl and A.~F. Vasseur.
\newblock{New perspectives in fluid dynamics: mathematical analysis of a model proposed by Howard Brenner.}
\newblock{\em New directions in mathematical fluid mechanics, 153–179, Adv. Math. Fluid Mech,} Birkhäuser Verlag, Basel, 2010.
\bibitem{FMT}  
U. S. Fjordholm, S. Mishra, and E. Tadmor. 
\newblock { On the computation of measure-valued solutions}. 
\newblock {\em Acta Numer.}, 25:567--679, 2016.  
  
  \bibitem{GKS}
   P. Gwiazda, O. Kreml, and A. \'Swierczewska-Gwiazda. 
   \newblock  { Dissipative measure-valued solutions for general conservation laws.}
   \newblock {\em Annales de l’Institut Henri Poincar\'e C, Analyse non lin\'eaire,} 37(3):683--707, 2020.
  
   \bibitem{GSW}
 P. Gwiazda, A. \'Swierczewska-Gwiazda, and E. Wiedemann. 
 \newblock {Weak-strong uniqueness for measure-valued solutions of some compressible fluid models. }
 \newblock {\em Nonlinearity}, 28(11):3873--3890, 2015.

 \bibitem{PAM}
B.~Haspot.
\newblock { From the highly compressible {N}avier--{S}tokes equations to fast
  diffusion and porous media equations, existence of global weak solution for
  the quasi-solutions.}
\newblock {\em J. Math. Fluid Mech.}, 18(2):243--291, 2016.
 
 \bibitem{HZ}
B.~Haspot and E.~Zatorska.
\newblock { From the highly compressible {N}avier-{S}tokes equations to the
  porous medium equation -- rate of convergence.}
\newblock {\em Discrete Contin. Dyn. Syst.}, 36(6):3107--3123, 2016.
	
	\bibitem{LSU}
O.~A. Lady{\v{z}}enskaja, V.~A. Solonnikov, and N.~N. Uralceva.
\newblock {\em Linear and quasilinear equations of parabolic type}.
\newblock Translated from the Russian by S. Smith. Translations of Mathematical
  Monographs, Vol. 23. American Mathematical Society, Providence, R.I., 1967.
  

\bibitem{LW}
M. J. Lighthill and J. B. Whitham. 
\newblock {On kinematic waves: I. Flow movement
in long rivers. II. A theory of traffic flow on long crowded roads.}
\newblock{\em Proc. Roy. Soc.}, A229:1749--1766, 1955.
  
\bibitem{PLL} 
P.-L. Lions. 
\it Mathematical Topics in Fluid Mechanics, Vol 2: Compressible Models. 
\rm Oxford University Press, \rm New York, 1998.
	
\bibitem{MNRR}
J. M\' alek, J. Ne\v cas, M. Rokyta, and M. R\r{u}\v zi\v cka. 
\newblock {\em Weak and measure-valued solutions to evolutionary PDEs.}
\newblock Volume 13 of Applied Mathematics and Mathematical Computation. Chapman and Hall, London, 1996.
	
\bibitem{MPI}	
S. Mueller
\newblock {\em	Variational  models for microstructure and phase transitions.} 
\newblock{ Calculus of Variation and Geometric Evolution Problem, Lecture Notes in Math., vol. 1713, Springer-Verlag, Berlin Heidelberg (1999)}
	
\bibitem{NS}
A.~Novotn{{\'y}} and I.~Stra{\v{s}}kraba.
\newblock {\em Introduction to the mathematical theory of compressible flow.}
  Volume~27 of {\em Oxford Lecture Series in Mathematics and its Applications}.
\newblock Oxford University Press, Oxford, 2004.
	

	




\bibitem{Payne}
H. J. Payne.
\newblock { FREFLO: A macroscopic simulation model of freeway traffic.}
\newblock {\em Transportation Research Record}, 722:68--75, 1979.


\bibitem{P}
P. Pedregal. 
\newblock{\em Parametrized measures and variational principles.}
\newblock{ Progress in Nonlinear Differential Equations and their Applications, 30.} Birkhäuser Verlag, Basel, 1997.

\bibitem{R}
P. I. Richards.
\newblock { Shock waves on the highway.}
\newblock {\em Operations Research}, 4:42--51, 1956.

\bibitem{tad4}
R.~Shvydkoy and E.~Tadmor.
\newblock { Eulerian dynamics with a commutator forcing {III}. {F}ractional diffusion of order $0<\alpha<1$.}
\newblock {\em Phys. D}, 376-377:131 -- 137, 2018.




\bibitem{VaYu}
A.~F. Vasseur and C.~Yu.
\newblock {Existence of global weak solutions for 3d degenerate compressible {N}avier--{S}tokes equations.}
\newblock {\em Invent. Math.}, 206(3)974, 2016.




\bibitem{W2018}
E. Wiedemann. 
\newblock  { Weak-strong uniqueness in fluid dynamics}. 
\newblock {\em London Math. Soc. Lecture Note Ser.,} pages 289–326. Cambridge Univ. Press, Cambridge, 2018.

\bibitem{Whitham}
G. B. Whitham.
\newblock {\em Linear and nonlinear waves.}
\newblock Wiley, New York, 1974.

\bibitem{Wo}
J. Wo\'znicki.
\newblock{Mv-strong uniqueness for density dependent, non-Newtonian, incompressible fluids.}
\newblock {\em arXiv:2101.02263, 2021.}

 \end{thebibliography}
\end{document}